\documentclass[12pt,reqno]{amsart}

\usepackage{latexsym}

\renewcommand{\Re}{{\operatorname{Re}}}
\renewcommand{\Im}{{\operatorname{Im}}}

\newcommand{\Morse}{{\operatorname{Morse}}}
\newcommand{\indx}{{\operatorname{index}}}
\newcommand{\diag}{{\operatorname{diag}}}
\newcommand{\sym}{{\operatorname{Sym}}}
\newcommand{\dual}{^{\scriptstyle *}}

\newcommand{\hol}{{\operatorname{Hol}}}
\newcommand{\crit}{{\operatorname {crit}}}
\newcommand{\K}{{\mathbf K}^\crit}
\renewcommand{\k}{{\mathbf k}^\crit}
\newcommand{\Hess}{{\operatorname {Hess}}}

\newcommand{\szego}{Szeg\"o }

\newcommand{\inv}{^{-1}}
\newcommand{\kahler}{K\"ahler }

\newcommand{\wt}{\widetilde}
\newcommand{\wh}{\widehat}
\newcommand{\PP}{{\mathbb P}}
\newcommand{\R}{{\mathbb R}}
\newcommand{\C}{{\mathbb C}}

\newcommand{\Z}{{\mathbb Z}}

\newcommand{\CP}{\C\PP}
\renewcommand{\d}{\partial}
\newcommand{\dbar}{\bar\partial}
\newcommand{\ddbar}{\partial\dbar}
\newcommand{\U}{{\rm U}}
\newcommand{\D}{{\mathbf D}}
\newcommand{\E}{{\mathbf E}\,}

\newcommand{\half}{{\frac{1}{2}}}
\newcommand{\vol}{{\operatorname{Vol}}}

\newcommand{\FS}{{{\operatorname{FS}}}}

\renewcommand{\phi}{\varphi}

\newcommand{\ccal}{\mathcal{C}}
\newcommand{\dcal}{\mathcal{D}}

\newcommand{\ical}{\mathcal{I}}
\newcommand{\kcal}{\mathcal{K}}
\newcommand{\lcal}{\mathcal{L}}
\newcommand{\ncal}{\mathcal{N}}

\newcommand{\ocal}{\mathcal{O}}
\newcommand{\scal}{\mathcal{S}}

\newcommand{\jcal}{\mathcal{J}}

\newcommand{\ga}{\gamma}

\newcommand{\La}{\Lambda}
\newcommand{\la}{\lambda}
\newcommand{\ep}{\varepsilon}
\newcommand{\de}{\delta}
\newcommand{\De}{\Delta}
\newcommand{\om}{\omega}

\newcommand{\CC}{{\mathcal C}}

\newcommand{\CL}{{\mathcal L}}

\newcommand{\CN}{{\mathcal N}}

\newcommand{\p}{{\partial}}

\newcommand{\GeV}{{\rm GeV}}
\newtheorem{maintheo}{{\sc Theorem}}
\newtheorem{maincor}[maintheo]{{\sc Corollary}}
\newtheorem{theo}{{\sc Theorem}}[section]

\newtheorem{lem}[theo]{{\sc Lemma}}

\newenvironment{rem}{\medskip\noindent{\it Remark:\/} }{\medskip}
\newtheorem{defin}[theo]{\sc Definition}

\title{Critical points and supersymmetric vacua}

\author{Michael R. Douglas}
\address{Rutgers, IHES and Caltech}
\email{mrd@physics.rutgers.edu}

\author{Bernard Shiffman}
\address{Department of Mathematics, Johns Hopkins University, Baltimore,
MD 21218, USA} \email{shiffman@math.jhu.edu}

\author{Steve Zelditch}
\address{Department of Mathematics, Johns Hopkins University, Baltimore,
MD 21218, USA} \email{zelditch@math.jhu.edu}

\thanks{Research partially supported by
DOE grant DE-FG02-96ER40959  (first author) and NSF grants
DMS-0100474 (second author) and DMS-0302518 (third author).}

\date{February 19, 2004; revised June 2, 2004}

\begin{document}

\begin{abstract} Supersymmetric vacua (`universes') of string/M theory may be
identified with certain critical points of a holomorphic section
(the `superpotential')  of a Hermitian holomorphic line bundle
over a complex manifold. An important physical problem is to
determine how many vacua there are and how they are distributed.
The present paper initiates the study of the statistics of
critical points $\nabla s = 0$ of Gaussian random holomorphic
sections with respect to a connection $\nabla$. Even the expected
number of critical points depends on the curvature of $\nabla$.
The principal results give formulas for the expected density and
number of critical points of Gaussian random sections relative to
$\nabla$ in a variety of settings. The results are particularly
concrete for Riemann surfaces. Analogous results on the density of
critical points of fixed Morse index are given.
\end{abstract}

\maketitle

\ \vspace{-.25in} \tableofcontents

\section{Introduction}

In a recent series of articles, one of the authors (M. R. Douglas)
has initiated a program to study the vacuum selection problem of
string/M theory from a statistical point of view \cite{Doug, AD}.
This is the problem of finding which (if any) solutions of the
theory describe our universe.  At present, we know very little
about the nature of this problem: how many solutions there are,
how many solutions we should expect to fit present-day
observations, the range of further predictions these solutions
make, and so forth.  But, we are reaching the point where such
questions can be studied systematically.

We briefly summarize the standard paradigm, called
compactification. String theory is formulated in ten space-time
dimensions.  To describe our universe, one considers a solution
which is a direct product of four dimensional Minkowski
space-time, with a compact six dimensional manifold (say a
Calabi-Yau manifold with metric) carrying additional structures (a
vector bundle and other discrete choices).  Given such a solution,
one can derive an ``effective field theory,'' in which the
continuous parameters of the solution (such as the complex
structure moduli of the Calabi-Yau) become fields in four
dimensional space-time.  One then looks for ``vacuum
configurations'' (or simply ``vacua'') of this field theory,
meaning stable time-independent classical solutions; each is a
possible candidate to describe our universe.  The physical
properties of a vacuum (such as masses of particles) are then
obtained by the analysis of small fluctuations around the
solution.

In the best studied (and perhaps most plausible) examples, the
effective field theory is an ``$N=1$ supergravity theory''
\cite{WB}. A primary datum in this theory is a {\it
superpotential} $W$, a holomorphic section of a line bundle $\CL
\to {\mathcal M}_{CY}$ over the moduli space of Calabi-Yau
manifolds.  A large class of vacuum configurations, the
supersymmetric vacua, are the critical points of this
superpotential.  Thus, part of the problem of counting vacua  is
to count critical points of a given holomorphic section.

Now at present there is no computation of an exact superpotential
for any string compactification.  Thus, to get insight into
general features of this problem, one might model the
superpotential as a random holomorphic section of $\CL$, much as
the potential in a disordered system is regarded as random, and
study the statistics of vacua (critical points) of a random
superpotential.

In setting up a statistical model of `random superpotentials', one
must decide which probability measure to put on the space of
candidates.  From many points of view, the simplest candidate is
to take Gaussian random holomorphic sections in $H^0({\mathcal
M}_{CY}, {\mathcal L})$ with respect to a certain covariance
kernel. This connects the string/M problem with the statistical
theory of holomorphic sections developed by the other two authors
in collaboration with P. Bleher in \cite{SZ, BSZ2, BSZ3}.

A further point in favor of this probability measure, is that it
can be used to get physical results.  Now there are known
superpotentials which become exact upon taking partial limits in
moduli space; for example taking the volume of the Calabi-Yau
metric large, at arbitrary complex structure.  A particularly
interesting class of these are the ``flux superpotentials''  which
are linear combinations from a computable basis of sections.  As
argued in \cite{AD}, an asymptotic estimate for the number of flux
vacua in a family of compactifications constructed by Giddings,
Kachru and Polchinski \cite{GKP} can be obtained from the expected
number of critical points in a Gaussian ensemble of
superpotentials.

We discuss the physical background of our problem in more depth in
section \ref{PB}, and for the remainder of the introduction,
concentrate on the mathematical results.

From the mathematical viewpoint, the statistical theory of
critical points of Gaussian random holomorphic sections depends on
the following objects:
\begin{itemize}

\item A choice of subspace $\scal \subset H^0(M, L)$  of holomorphic sections of a  holomorphic line
bundle $L \to M$. We assume $\scal$ to be finite dimensional, but
our methods and results extend easily to infinite dimensional
spaces of sections.

\item A choice of  Gaussian measure $\gamma$ on $\scal$, or equivalently
an inner product $\langle, \rangle$ on $\scal$.

\item A choice of Hermitian metric $h$ on $L$.
This gives rise to the Chern connection $\nabla$ on $L$, which is
of type $(1,0)$ with curvature of type $(1,1)$.

\end{itemize}

The physical application requires a rather general framework of
$(\scal, \gamma ,\nabla)$. Our purpose is to study the
distribution and number of critical points \begin{equation}
\label{CPE} \nabla s (z) = 0,\;\;\; s \in \scal \end{equation}
 of a random section $s \in \scal$ with respect to $\gamma$.
 Since critical points are  zeros of random sections $\nabla s
 \in \nabla \scal \subset \ccal^\infty(T^*M \otimes L)$,  we are able to use the  previous work \cite{SZ, BSZ2,
BSZ3} on zeros of random $\ccal^\infty$ sections of complex vector
bundles to set up the statistical theory. It is important to
observe
 that the critical point  equation (\ref{CPE}) is not holomorphic,
 and therefore the much simpler statistical theory of zeros of
 holomorphic sections in \cite{SZ, BSZ2, BSZ3}  does not apply. On the other hand,  the
 holomorphicity of the original sections permeates the calculations,
  and gives the statistics of their critical points some special
features which do not hold in the general $\ccal^\infty$ case.  The
purpose of this article is to develop a self-contained theory of
critical points of Gaussian random holomorphic sections which
makes use of these special features. It is not necessary to know
the results of \cite{ BSZ2, BSZ3} to read this paper. At the risk
of being repetitious, we have tried to make this article more
accessible to physicists as well as mathematicians by sometimes
giving two proofs of the same assertion or two explanations of a
key idea in both mathematical and physical language.

The key object of interest in this article is   the {\it expected
distribution of critical points}  of a
Gaussian random holomorphic  section $s \in \scal \subset H^0(M,
L)$. The inner product on $\scal$  determines a Gaussian
probability measure $\gamma$ (see Definition \ref{GM}). As
discussed in \S \ref{CPHS}, the  definition $\nabla s = 0$  of
critical point depends on a choice of connection, which  we always
choose  to be the  Chern (Hermitian) connection associated to a
Hermitian metric $h$ on $L$. For almost any section $s \in H^0(M,
L)$, the set of its critical points $Crit^{\nabla}(s)$ is discrete
and we define the distribution of critical points of $s$ to be the
(un-normalized)  measure
\begin{equation}\label{CSH}  C_s^{\nabla} := \sum_{z\in Crit^{\nabla}(s)}
\delta_{z},
\end{equation}
where $\delta_{z}$ is the Dirac point mass at $z$. We let
\begin{equation}\label{DENSITY}\K_{\scal,\ga,\nabla}=
\E_{\gamma} C_s^{\nabla}\end{equation} denote the expected
distribution of critical points, i.e.\ the average of the measures
$C_s^{\nabla}$ with respect to $\gamma_{\scal}$. If $M$ is
compact, the expected total number of critical points is then
given by
\begin{equation} \label{NUMBER0}
{\mathcal N}^{\crit}(\scal,\ga,\nabla) =
\K_{\scal,\ga,\nabla}(M)\;.
\end{equation}

Our first result is a formula for the expected critical point
distribution $\K_{\scal,\ga,\nabla}$, which is valid for any
subspace $\scal \subset H^0(M, L)$ of holomorphic sections of any
holomorphic line bundle over any complex  manifold (possibly
non-compact and/or incomplete). We shall assume that the space
$\scal$ of sections satisfies a technical condition, the {\it
2-jet spanning property\/}, which says that all possible values
and derivatives of order $\le 2$ are attained by the global
sections $s \in \scal$ at every point of $M$ (see Definition
\ref{2JET}). In particular, the formula for $\K_{\scal,
\gamma,\nabla}$ applies to the physically relevant case in
string/M theory where $L \to M$ is a negative line bundle over the
moduli space of Calabi-Yau manifolds (an incomplete, non-compact
\kahler manifold), where $\scal$ is a special subspace of sections
given by periods of the Calabi-Yau form, and where $\gamma$ is
induced by   a rather subtle inner product coming from the
intersection form on cycles.

To state the result, it is most convenient to introduce a local
frame (non-vanishing holomorphic section) $e_L$ for $L$  and local
coordinates $(z_1, \dots, z_m)$ on an open set $U\subset M$, and
to write
\begin{equation}\label{denfn}  \K_{\scal, \gamma, \nabla}=
\k_{\scal, \gamma, \nabla}\,dz\;,\end{equation} where
$dz=\prod_{j=1}^m\left(\frac i2 dz_j\wedge d\bar z_j\right)^m$ is
Lebesgue measure with respect to these coordinates. (Note that
$\k_{\scal, \gamma, \nabla}$ depends on the coordinates.) We
denote the curvature $(1,1)$-form of $\nabla$ in these coordinates
by $\Theta = \sum_{j,k = 1}^m \Theta_{jk} dz^j \wedge d\bar{z}^k$
and refer to the $m \times m$ matrix $\big(\Theta_{jk}\big)$ as
the curvature matrix of $\nabla$. We also denote by $\sym(m, \C)$
the space of complex $m \times m$ symmetric matrices. It is a
Hermitian vector space with inner product $\langle A, B \rangle =
Tr A B^*$. We also consider the Hermitian orthogonal sum $\sym(m,
\C) \oplus \C$ with the standard Hermitian inner product on $\C$.
These induce a natural volume form on $\sym  (m, \C) \oplus \C$.
Finally, the {\it covariance kernel} or two-point function
$\Pi_{\scal}(z,w)$ of the Gaussian measure $\gamma$ is defined in
Definition \ref{COVKER}.

\begin{maintheo}\label{KNcrit1} Let $(\scal,
\gamma, \nabla)$ denote a  finite-dimensional  subspace $\scal
\subset H^0(M, L)$ of holomorphic sections of a holomorphic line
bundle $L \to M$ with a Chern connection $\nabla$ on an
$m$-dimensional complex  manifold, together with a Gaussian
measure $\gamma$ on $\scal$. Assume that $\scal$ satisfies  the
$2$-jet spanning property.  Given local coordinates $z = (z_1,
\dots, z_m)$ and a local frame $e_L$ for $L$, there exist
positive-definite Hermitian matrices
$$A(z): \C^m \to \C^m\;, \quad
\Lambda(z): \sym  (m, \C) \oplus \C \to \sym (m, \C)\oplus
\C\;,$$ depending only on  $z$, $\nabla$ and
  $\Pi_{\scal}$ (cf.\
\eqref{AN1}--\eqref{CN1}) such that the expected density  of
critical points with respect to Lebesgue measure $dz$ is given by
\begin{eqnarray*}\k_{\scal,\gamma, \nabla} (z)  &=&
\frac{1}{\pi^{{m+2\choose 2}}\det A(z) \det\La(z)}\\ &&\quad\times
\int_{\C}\int_{\sym (m, \C)}  \left|\det\begin{pmatrix}
H' &x\,\Theta(z)\\
\bar x\,\bar\Theta(z)&\bar H'\end{pmatrix} \right|\; e^{
-{\left\langle \La(z)^{-1}(H'\oplus x), \,H'\oplus x
\right\rangle}}\,dH'\,dx\,,
\end{eqnarray*}
where $\Theta(z)$ is the curvature matrix of $\nabla$ in the
coordinates $(z_1, \dots, z_m)$.
\end{maintheo}
The matrix in the formula is the complex Hessian of $s$ discussed
in \S \ref{HESS1}.

In order to give the simplest expressions for the matrices $A(z)$
and $\La(z)$ , we let $e_L$ be an {\it adapted\/} local frame at a
point $z_0\in M$; i.e., $e_L$ has the property that the pure
holomorphic derivatives (of order $\le 2$) of the local connection
form for $\nabla$ vanish at $z_0$ (see Definition \ref{ADAPTED}).
 We then let $F_{\scal}(z,w)$ be the
local expression for $\Pi_{\scal}(z,w)$ in the frame $e_L$ (see
Definition \ref{xexpect}). Then  \begin{equation} A(z_0)=
\left(\left. \frac{\partial^2}{\partial z_j \partial \bar{w}_{j'}}
 F_{\scal} (z,w)\right|_{(z,w)=(z_0,z_0)}\right) \label{AN1}\end{equation}
and
\begin{equation}
\La(z_0)=C(z_0)-B(z_0)^*A(z_0)^{-1}B(z_0) \;,\label{Lambda}
\end{equation}
where
\begin{eqnarray}
B(z_0)&=&\left.\left[\left(  \frac{\partial^3}{\partial z_j \partial
\bar{w}_{q'} \partial \bar{w}_{j'}} F_{\scal} (z,w)\right)
\quad \left(  \frac{\partial}{\partial z_j}
F_{\scal} (z,w)\right)
\right]\right|_{(z,w)=(z_0,z_0)} \,,\label{BN1}\\[8pt]
C(z_0)&=&\left.\left[
\begin{array}{cc}\left(
\frac{\partial^4}{\partial z_q \partial z_j \partial
\bar{w}_{q'}\partial \bar{w}_{j'}} F_{\scal} (z,w) \right) &
\left(  \frac{\partial^2}{\partial z_j
\partial z_q}F_{\scal} (z,w)
\right)
\\[8pt]
\left(  \frac{\partial^2}{\partial \bar{w}_{q'}\partial
\bar{w}_{j'}} F_{\scal} (z,w)\right) & F_{\scal}(z,z)
\end{array}\right]\right|_{(z,w)=(z_0,z_0)} \,,\label{CN1}\\ &&
\qquad\qquad \qquad 1\le j\le m\,, 1\le j\le q\le m\,, 1\le j'\le q'\le
m\,.\nonumber
\end{eqnarray}

In the above, $A,B,C$ are $m\times m,\,m\times n,\, n\times n$
matrices, respectively, where $n=\half(m^2+m+2)$. We also provide
formulas for these matrices in a non-adapted local frame (cf.\
(\ref{AN})--(\ref{CN})), which are useful when one studies
variations of the critical point distribution with respect to
$\nabla$.

Of course, the expected distribution of zeros
$\kcal^\crit_{\scal,\gamma, \nabla} (z) dz$ is independent of the
choice of frame and coordinates.  One can write the formula in an
invariant form by interpreting
$$A(z):T^{*1,0}_{M,z}\otimes L_z \to T^{*1,0}_{M,z}\otimes L_z\;, \quad
\Lambda(z):(S^2T^{*1,0}_{M,z}\oplus \C)\otimes L_z \to
(S^2T^{*1,0}_{M,z}\oplus \C)\otimes L_z\;$$ as positive-definite
Hermitian operators, where $S^2T^{*1,0}_{M,z}\subset
T^{*1,0}_{M,z}\otimes T^{*1,0}_{M,z}$ denotes the symmetric
product, $\Theta(z) \in T^{*1,0}_{M,z}\otimes T^{*0,1}_{M,z}$ is
the curvature operator of $\nabla$, and the determinant in the
integral is an element of $(\det T^*_{M,z})^2\otimes L^m_z$.
However, we find the local expressions to be more useful.

Since the integral contains an absolute value, it is difficult to
evaluate the density explicitly when the dimension is greater than
1, or even to analyze its dependence on $\nabla, \gamma$. In
particular, one cannot simplify it with Wick's formula.

A special case of geometric interest is  where the inner product
$\langle, \rangle$ and Gaussian measure on $\scal$ are induced by  a volume
form  $dV$ on $M$ and  the same  Hermitian metric $h$ which determines
$\nabla$, namely
\begin{equation} \label{HIP} \langle s_1, s_2 \rangle = \int _M h_z\big(s_1(z),s_2(z)\big)
\,dV(z)\,. \end{equation} The covariance
kernel is then the
\szego kernel of $\scal$, i.e. the orthogonal projection
\begin{equation}\label{szegoproj}
(\Pi_{\scal, h, V}\, s)\,(z) = \int_M h_w\big(s(w),\Pi_{\scal, h,
V}(z, w)\big) dV(w)\,,\qquad \Pi_{\scal,h, V}(z,w)\in L_z\otimes
\bar L_w\;.
\end{equation}
We refer to this Gaussian measure as the {\it Hermitian Gaussian
measure} on $\scal$. In this case, every object in the density of
critical points is determined by the metric and volume form, and
we have a direct relation between the expected number and
distribution of critical points and the metric.

 The
simplest metric situation   is that of positive   line bundles. In
this case, we assume that $\om=\frac{i}{2}\Theta_h$ so that
$c_1(L) = [\frac 1\pi \omega]$,
 where the brackets denote the cohomology
class. More precisely, the \kahler form is given by
$$\om=\frac{i}{2}\Theta_h= \frac i2 \ddbar K, \qquad K =-\log|e_L|^2_h. $$
The volume form is then assumed to be $$dV = \frac{\omega^m}{m!}$$
(and thus the total volume of $M$ is
$\frac{\pi^m}{m!}c_1(L)^m$\,).  If $L$ is a negative line bundle
(on a noncompact manifold $M$), we instead choose the \kahler form
$\om=-\frac{i}{2}\Theta_h$. We shall write
\begin{equation} \K_{\scal,\ga,\nabla} =
\kcal^\crit_{\scal,h}\,dV\;,\;\;\;\; {\mathcal N}^{\crit}(\scal,
\gamma,  \nabla ) = {\mathcal N}^{\crit}(\scal, h),
\end{equation} where $\ga$ is the Hermitian Gaussian measure
described above. (Note that $\kcal^\crit_{\scal,h}$ denotes the
density with respect to the volume $dV$, while
$\k_{\scal,\ga,\nabla}$ is the density with respect to Lebesgue
measure in local coordinates.) As a consequence of Theorem
\ref{KNcrit1}, we obtain the following integral formula for the
critical point density $\kcal^\crit_{\scal,h}$ in these cases:

\begin{maincor}\label{KNcrit2} Let
$(L,h)\to M,\ (\scal,  \gamma,\nabla)$ be as in Theorem
\ref{KNcrit1}. Further assume that the curvature form of
$\Theta_h$ is  either positive or negative and that $\gamma$ is
the Hermitian Gaussian measure \eqref{HIP} induced by $h$  and by
the volume form $dV= \frac 1{m!} \left(\pm \frac
i2\Theta_h\right)^m$. Then  the expected density of critical
points relative to $dV$ is
 given by
$$\kcal^\crit_{ \scal,h} (z) =
\frac{\pi^{-{m+2\choose 2}}}{\det A(z) \det\La(z)} \int_{\sym(m,\C)
\times \C} \left|\det(H'H'{}^*-|x|^2I)\right|e^{ -{\langle
\La(z)^{-1}(H', x),(H', x) \rangle}}\,dH'\,dx\,. $$
where $A(z), \Lambda(z)$ are positive Hermitian matrices
 (depending on $h$ and $z$) given by
\eqref{AN1}--\eqref{CN1}.
\end{maincor}
As above,  $H' \in \sym(m, \C)$ is a complex symmetric matrix, and
the matrix $\Lambda$ is a Hermitian operator on the complex vector
space $\sym(m, \C) \times \C$. The only point of the corollary  is
that we can identify $\Theta = I$ in normal coordinates and
simplify the determinant.

As mentioned above, the physically relevant case is that of a
negative line bundle over an incomplete \kahler manifold (the
moduli space of Calabi-Yau metrics on a $3$-fold $M$). As
mentioned above, the relevant  Gaussian measure there is not the
Hermitian one. However, the same  formula holds in that case since
the curvature form equals  $- I$ in local coordinates adapted to
the Weil-Petersson volume form.

In dimension one, we obtain the following explicit formula for the
expected density of critical points in terms of the eigenvalues of
$\La Q_r$, where
$$Q_r=\begin{pmatrix}1&0\\ 0& -r^2\end{pmatrix}\;,$$ and $r=\frac
i2
\Theta_h/dV$:

\begin{maintheo} \label{crit1} Let $(L, h) \to M$ be a  Hermitian
line bundle on a (possibly non-compact) Riemann surface $M$ with
area form $dV$. Let $\scal$ be a finite-dimensional subspace of
$H^0(M,L)$ with the 2-jet spanning property, and let $\ga$ be the
induced Hermitian Gaussian measure. Let
$\mu_1=\mu_1(z),\,\mu_2=\mu_2(z)$ denote the eigenvalues of
$\La(z) Q_r$, where $r=r(z)=\frac i2 \Theta_h/dV$. Then
$\K_{\scal,\ga,\nabla}=\kcal_{\scal,h,V}\,dV$, where
$$ \kcal^\crit _{\scal,h,V} = \frac{1}{\pi A}\; \frac{\mu_1^2 +
\mu_2^2 }{| \mu_1| + | \mu_2|}=  \frac{1}{\pi A}\; \frac{Tr\,
\La^2}{Tr | \La^\half Q_r\La^\half|}\;,
$$ where $A,\La$ are given by
\eqref{AN1}--\eqref{CN1}.
\end{maintheo}

We define the {\it topological index\/} of a section $s$  at a
critical point $z_0$ to be the index of the vector field $\nabla
s$ at $z_0$ (where $\nabla s$ vanishes). Critical points of a
section $s$ in dimension one are (almost surely) of topological
index $\pm 1$. (If the connection were flat, then $\nabla s$ would
be holomorphic and the topological indices would all be positive.)
The critical points of $s$ of index $1$ are the saddle points of
$\log |s|_h$ (or equivalently, of $|s|^2_h$), while those of
topological index $-1$ are local maxima of $\log|s|_h$ in the case
where $L$ is positive, and are local minima of $|s|^2_h$ if $L$ is
negative. (If $L$ is negative, the length $|s|_h$ cannot have
local maxima; if $L$ is positive, the only local minima of $|s|_h$
are where $s$ vanishes.)  Thus, in dimension 1, topological index
$1$ corresponds to $\log|s|_h$ having Morse index 1, while
topological index $-1$ corresponds to Morse index 2 if $L$ is
positive.  In fact, in all dimensions, the critical points of a
section $s$ are the critical points of $\log|s|_h$, and for
positive line bundles $L$, we have
\begin{equation}\label{2indices} \mbox {index}_{z_0}(\nabla s) =
(-1)^{m+\Morse \; \indx_{z_0} (\log |s|)}\;,\end{equation}  at
(nondegenerate) critical points $z_0$ (see Lemma \ref{index}).

From the proof of Theorem \ref{crit1} we obtain:

\begin{maincor} \label{indexdensity} Let $(L, h) \to (M,dV),\
\mu_1,\mu_2$ be as in Theorem \ref{crit1}. Then:

\begin{itemize}

\item  The  expected density of critical
points of topological index $1$ (where $|s|^2_h$ has a saddle
point) is given by
$$\kcal^\crit_+(z)=
\frac{1}{\pi A(z)}\; \frac{\mu_1^2 }{| \mu_1| + | \mu_2|} \;,$$

\item The expected density of critical points of topological
index $-1$ (where $|s|^2_h$ has a local maximum) is
$$\kcal^\crit_-(z)=
\frac{1}{\pi A(z)}\;\frac{\mu_2^2  }{| \mu_1| + | \mu_2|} \;.$$

\item
Hence, the index density is given by
\begin{equation*}\kcal^\crit_{\rm index}: = \kcal^\crit_+(z) -
\kcal^\crit_-(z) = \frac{1}{\pi A(z)}\;
(\mu_1+\mu_2) = \frac{1}{\pi A(z)}\;Tr[\La(z)Q]
\;.\end{equation*}
\end{itemize}
\end{maincor}

The simplest case of Theorem \ref{crit1} and Corollary \ref{indexdensity}
is that of sections of powers $\ocal(N)$ of the hyperplane line bundle
$\ocal(1) \to
\CP^1$, i.e. of homogeneous polynomials of degree $N$ in the
$SU(2)$ ensemble, where we can give exact formulas:

\begin{maincor}\label{exactP1} The expected numbers
$\ncal^\crit_{N,+}$ and $\ncal^\crit_{N,-}$ of critical points of
topological index $1$ and   $-1$, respectively, of a random
section $s_N\in H^0(\CP^1,\ocal(N))$ (endowed with the Hermitian
Gaussian measure induced from the Fubini-Study metrics $h^N$ on
$\ocal(N)$ and $\om_\FS$ on $\CP^1$) are given by
\begin{eqnarray*}\ncal^\crit_{N,+}&=&
\frac{4(N-1)^2}{3N-2}\  =\ \frac 43\, N - \frac{16}{9} + \frac
4{27}\, N^{-1}\cdots \quad (\mbox{number of saddle points of }
|s|^2_h)\;,
\\\ncal^\crit_{N,-}&=&\frac{N^2}{3N-2}\  =\ \frac 13\, N +
\frac{2}{9} + \frac 4{27}\, N^{-1} \cdots \quad (\mbox{number of
local maxima of } |s|^2_h)\;,
\end{eqnarray*} and thus the expected total number of critical
points is given by
$${\mathcal N}^\crit_N(\CP^1)=\frac{5N^2-8N+4}{3N-2}  = \frac 53\, N -
\frac{14}{9} + \frac 8{27}\, N^{-1} \cdots\;.$$
\end{maincor}

\medskip  It follows  that the average number of critical points
of a polynomial $p(z)$ of degree $N>1$ in the $SU(2)$  ensemble is
greater for all $N$  than the almost sure number of critical
points (= $N - 1$) in the classical sense of $p'(z) = 0$. This is
not surprising, since in the former case, sections may have
critical points of index $-1$ and in the latter case there are no
critical points of index $-1$, while $\ncal^\crit_{N,+} -
\ncal^\crit_{N,-} = c_1\left(\ocal(N)\otimes K_{\CP^1}\right)=N-2$.
(The number of critical points in the latter case is $N-1$ instead
of $N-2$, since $p'(z)\,dz$ almost surely has a pole of order 1 at
$\infty$.) In an asymptotic sense, there are $\frac{5}{3}$ as many
critical points in the metric sense.

In higher dimensions, the integral in Corollary \ref{KNcrit2} is more complicated
to evaluate, and the density does not have a simple formulation in
terms of eigenvalues as in Theorem \ref{crit1}.  In our subsequent paper \cite{DSZ},
we derive the following alternate formula for the
 expected density of critical points:
\begin{equation*}
\kcal^\crit_{\scal,h}(z) =  \frac{c_m}{\det
A(z)}\lim_{\ep,\ep' \to 0^+}  \int_{\R^m }\int_{\R^m}
\int_{\U(m)}
\frac{\Delta(\xi)\,
\Delta(\lambda)\; |\prod_j \la_j| \,e^{i\langle \xi, \lambda
\rangle}   e^{- \epsilon |\xi|^2 -\epsilon'
|\lambda|^2}}{\det\left[ i\wh D(\xi)\rho(g)\La(z)\rho(g)^* +
I\right] }\, dg\, d \xi\, d\lambda,
\end{equation*}
 where
\begin{itemize}
\item $\displaystyle c_m=\frac{(-i)^{m(m-1)/2}}{2^m\,\pi^{2m}\,
\prod_{j=1}^mj!}\;,$
\item $\wh D(\xi)$ is
the Hermitian operator on $\sym(m,\C)\oplus \C$ given by
$$ \wh D(\xi) \big((H_{jk}),x\big) =
\left( \left(\frac{\xi_j+\xi_k}{2}\,H_{jk}\right),\
-\left(\textstyle\sum_{q=1}^m \xi_q\right) x \right)\;,$$ \item
 $\rho$ is the representation  of $\U(m)$ on
$\sym(m,\C)\oplus \C$ given by
$$\rho(g)(H, x) = (gHg^t,x)\;.$$
\end{itemize}
We use this formula in \cite{DSZ} to compute (with the assistance
of Maple) the expected numbers
$\ncal^\crit_N(\CP^m)$ of critical points of sections of the $N$-th power
$\ocal(N)$ of the hyperplane section bundle on $\CP^m$, for $m\le 4
$. For example,   for the projective plane, we have
\begin{equation}\label{CP2}\ncal^\crit_N(\CP^2)=
{\frac {59\,{N}^{5}-231\,{N}^{4}+375\,{N}^{3}-310\,{N}^{2}+132\,N-24}{
 \left( 3\,N-2 \right) ^{3}}} \sim \frac{59}{27}\,N^2\;.\end{equation}
This expected number is with respect to the $SU(3)$-invariant geometry on
$\ocal(N)\to\CP^2$.  We conjecture that \eqref{CP2}  is the minimum expected number
of critical points over all connections on $\ocal(N)\to\CP^2$ with positive
curvature, and we show in
\cite{DSZ} that this is the case in an asymptotic sense.

In our sequel \cite{DSZ}, we also  obtain asymptotic results on the expected
numbers of critical points of $f=|s|^2_h$ of each possible Morse
index for powers $L^N$ of a positive line bundle $(L,h)$ over any
compact complex manifold.  In particular, for a positive line
bundle $(L,h)$ over a compact Riemann surface $C$ endowed with the
\kahler form $\om_h=\frac i2 \Theta_h$, we prove that
\begin{eqnarray}\ncal^\crit_{N,+} &=& \frac 43 \,c_1(L)\,
N + \frac 89\, ( 2g-2) +
\left(\frac 1{27\pi}\,
\int_C \rho^2 \om_h\right)N\inv +O(N^{-2})\;,\label{critRS1}\\
\ncal^\crit_{N,-} &=& \frac 13 \,c_1(L)\, N - \frac 19\, ( 2g-2) +
\left(\frac 1{27\pi}\, \int_C \rho^2 \om_h\right)N\inv
+O(N^{-2})\;,\label{critRS2}\end{eqnarray} where $g$ denotes the
genus of $C$ and $\rho$ is the Gaussian curvature of $(C,\om_h)$.
Thus, the expected number of local maxima of $|s|_h$ (on Riemann
surfaces of any genus) is $\sim \frac{1}{3}c_1(L)N\sim \frac 13
\dim H^0(M,L^N)$. It would be interesting to find a heuristic
reason for the factor $\frac 13$.

It is well known that on complex manifolds $M$ of dimension $m$,
critical points of a section $s\in H^0(M,L)$ are critical points
of the function $\log|s|_h$ (and conversely) and these have Morse
index $\ge m$ when $L$ has positive curvature (see \cite{B}).
(Recall that the Morse index of a nondegenerate critical point of
a real-valued function is  the number of  negative eigenvalues of
its Hessian matrix.) In this case, we have a density formula for
critical points of any Morse index:

\begin{maintheo} \label{Morse} Let $(L, h) \to M$ be a positive
holomorphic line bundle over a complex manifold $M$ with volume
form $dV=\frac 1{m!}(\frac i2\Theta_h)^m$. Suppose that $H^0(M,
L)$ contains a finite-dimensional subspace $\scal$  with the
$2$-jet spanning property, and let $\ga$ be the Hermitian Gaussian
measure on $\scal$. Then the expected density  with respect to
$dV$ of critical points  of $\log|s|_h$ of Morse index $q$ is
given by
$$\kcal^\crit_{\scal,h,q} (z) =
\frac{\pi^{-{m+2\choose 2}}}{\det A(z) \det\La(z)} \int_{{\bf S}_{m,q-m}}
\left|\det(SS{}^*-|x|^2I)\right|e^{ -{\langle
\La(z)^{-1}(S, x),(S, x) \rangle}}\,dS\,dx\,. $$ where
$${\bf S}_{m,k}=\{S\in \sym(m,\C)\times \C: \mbox{\rm index}(SS^*-|x|^2I) =k\}\;.$$
\end{maintheo}

 This article is just the first  in a series and leaves many
 issues unexplored. For simplicity let us assume that $\scal = H^0(M, L)$ and drop it
 from the notation. First, it  would be
 interesting to explore the  dependence of the density and expected number
 ${\mathcal N}^{\crit}(\nabla, \gamma)$  of
critical points on the connection $\nabla$ and the Gaussian
measure $\gamma$. In \S \ref{s-dependence} we give a simple proof
that ${\mathcal N}^{crit}(\nabla, \gamma)$ is non-constant in
$(\nabla, \gamma)$. The number ${\mathcal N}^{\crit}(\nabla,
\gamma)$ might be viewed as defining a configurational entropy for
statistics of vacua. It is bounded below by the Euler
characteristic $c_m(L \otimes K_M)$, but (as the argument in \S
\ref{s-dependence} suggests) is probably not bounded above. It
would be interesting to prove this, and to analyze how ${\mathcal
N}^{\crit}(\nabla, \gamma)$ depends on the  choice of
 $\nabla, \gamma$? Does there exist a
smooth $\nabla$ with a minimal average number of critical points?
(Clearly, the meromorphic $\nabla$ gives the minimal number for
$\ocal(N) \to \CP^m$, but the number jumps as one moves from a
smooth to a meromorphic connection.)

These problems become  purely geometrical when the Gaussian
measure is  Hermitian, i.e.  $\gamma = \gamma_{\scal, h, V}$. In
this case,  does there exist a metric for which ${\mathcal
N}^{\crit}(h)$ is minimal? Is it unique?  The formulas above
express $ {\mathcal N}^{\crit}(h)$ in terms of the metric \szego
kernel and are not explicit in terms of the geometry of $(L, h)$.
How does  the curvature of $h$ influence  ${\mathcal
N}^{\crit}(h)$. Is it bounded in a set of metrics with curvature
bounds? Do positively curved metrics on $L$ have larger ${\mathcal
N}^{\crit}(h)$ than signed curvature ones, as suggested by the
uniform distribution result? Can one link
 the expected number of critical points for meromorphic
 connections with that for smooth connections?

 In
subsequent work \cite{DSZ}, we will analyze the asymptotics of the
density and number ${\mathcal N}^{\crit}(h^N)$ of critical points
for powers $(L^N, h^N) \to (M, \omega)$ of  a positive Hermitian
line bundle $(L, h) \to (M, \omega)$.  We will show that, as in
Theorem \ref{exactP1}, the
 density $\kcal_{N, h}^{\crit}$ of critical points has a complete asymptotic
 expansion in $N$, whose leading coefficient  is a universal
 constant  with respect to the curvature volume form $\frac{\omega^m}{m!}$.
 Thus, critical points become uniformly distributed with respect
to the curvature volume form. Furthermore, we will analyze the
asymptotic dependence of the expected number of critical points
${\mathcal N}(h^N)$ on the metric. But this asymptotic study does
not seem to answer the above questions on a fixed positive line
bundle.

 Finally, the motivating problem is that of statistics of vacua in
 string/M theory. In this case, the  number  ${\mathcal
N}^{\crit}(\scal, \nabla, \gamma)$ depends  on the particular
choice of the ample subspace $\scal$  of periods and on the
special choice of $\gamma$ coming from the intersection form.
Moreover, the Gaussian measure is only an approximation to the
discrete probability space of periods. Ultimately, we would like
to understand the above questions in this setting.

\section{Physical Background}\label{PB}

In this section, we give precise definitions for the physical
theories we study.  As stated in the introduction, we do not
discuss string or M theory directly, but rather assume that a
given string or M theory compactification corresponds to an
``effective $\CN=1$ supergravity theory,'' in a way we sketch in
an example below.

The standard references for supergravity and other field theories
with ``$\CN=1$ supersymmetry'' are \cite{WB,Weinberg}, and nice
treatments of supersymmetry for mathematicians are
\cite{Freed,IAS}.  Field theories are usually defined by
specifying an action functional, which is written in terms of
fields which are sections of various spinor and tensor bundles
over $\R^{D,1}$, taking values in a configuration space $M$ and
its associated bundles.

For present purposes, the basic data specifying a supergravity
theory $T$ is a triple $(M,K,W)$, where

\begin{itemize}

\item  $M$ is the ``configuration space,'' a complex K\"ahler
manifold.  We will typically denote its dimension as $d$, and
local complex coordinates as $z^i$.  We will also refer to these
coordinates as ``fields.''

\item
$K$ is the K\"ahler potential, determining the metric on $M$.

\item  $W$ is the superpotential, a holomorphic section (possibly
with singularities) of the associated line bundle $\CL$ with
$c_1(\CL)=-\frac 1\pi \omega$, where the K\"ahler form
$\omega=\frac i2 \ddbar K$.

\end{itemize}

Such an associated line bundle carries a natural holomorphic
connection whose curvature is the K\"ahler form. It is the
connection which preserves the Hermitian metric on the fibers
$$
||W||^2 = e^K |W|^2 ;
$$
explicitly, the covariant derivative of a section $W$ is
\begin{equation}\label{eq:covder}
D_i W \equiv \p_i W + (\p_i K) W ,
\end{equation}
while the curvature is
$$
F=[D,\bar \p] = -\omega .
$$
The importance of this structure was first emphasized in
\cite{BaggerWitten}.

From this data, one can construct the {\it scalar potential} $V$.
It is the following function on $M$ (\cite{WB} formula 21.22; p.
169):
\begin{equation}\label{eq:defV}
V = e^K \left( g^{i\bar j} (D_i W) (\bar D_{\bar j} W^*) - 3 |W|^2
\right)
\end{equation}
where $W^*(\bar z)$ is the complex conjugate section and $ \bar
D_{\bar j} W^* = \bar \p_{\bar j} W^* + (\bar \p_{\bar j} K) W^* .
$

The basic physics of the scalar potential is the following.  While
the fields $z^i$ are functions on four-dimensional space-time, in
a state of minimum energy (a ground state or vacuum) they will
take constant values, at which the potential energy function
$V(z)$ is a local minimum.  In a given compactification and its
corresponding supergravity theory $T$, there could be one, several
or no such minima.  The case of multiple minima is rather
analogous to the familiar phenomenon of ``phases of matter'' such
as solid, liquid, and gas, which have different expectation values
for position-independent ``fields'' such as the local density,
pressure, and so forth, which could be determined by minimizing a
free energy.  All physical predictions depend on the choice of
minimum, and a first step to understanding the consequences of
this is to know how many minima there are.

Whereas in general, the scalar potential in a field theory can be
an arbitrary real function. in supergravity it must take the form
(\ref{eq:defV}), so this is a key formula in the physics of
supersymmetry.  Some of its features admit a more conceptual
explanation.  For example, the important fact that it is
sesquilinear in $W$ with signature $(n,1)$, and is thus not
positive definite, is the expected generalization of the familiar
statement that a supersymmetric Hamiltonian is a sum of squares,
to a theory containing gravity.

We define a {\it vacuum} to be a critical point $p\in\CC$ of $V$.
The vacua are further distinguished as follows:

\begin{itemize}

\item A {\it supersymmetric} vacuum is one in which the covariant
gradient $D_i W=0$.  This can easily be seen to imply $V'=0$, but
the converse is not true.

\item A {\it non-supersymmetric} vacuum is a critical point
$V'=0$, at which $D_i W \ne 0$. The norm of the gradient,
\begin{equation}\label{eq:susyscale}
M_{susy}^4 \equiv e^K g^{i\bar j} D_i W D_{\bar j} W^* ,
\end{equation}
is then referred to as the {\it scale of supersymmetry breaking}.

\item The value of $V$ at a critical point is the {\it cosmological
constant} $\Lambda$ of that vacuum.  These are divided into
$\Lambda=0$, the Minkowski vacua, $\Lambda>0$, the {\it de Sitter}
(or dS) vacua, and $\Lambda<0$, the {\it Anti-de Sitter} (or AdS)
vacua. It is easy to see that supersymmetric vacua can only be
Minkowski or AdS.  The {\it Minkowski vacua} are simultaneous
solutions of $D_i W=W=0$; in this case $D_i W = \p_i W$ and the
existence of such vacua is independent of the K\"ahler potential.
On the other hand, this is an overdetermined set of equations, so
generic superpotentials do not have  supersymmetric Minkowski
vacua.

\end{itemize}

For our purposes, a {\it metastable} vacuum will be one for which
the Hessian $V''$ is non-negative definite.  Physically, this is
required so that small fluctuations of the fields will not grow
exponentially.\footnote{ This is evident in a Minkowski vacuum,
but not literally true in a supersymmetric AdS vacuum; however
this condition is still interesting in the latter context as it is
the condition for stability when supersymmetry is broken by other
(D term) effects.} We use the term metastable rather than stable,
as such vacua have other potential instabilities (tunnelling)
which we mention below.

We finally make a few comments about units.  As with general
relativity, in supergravity it is natural to work in ``Planck
units,'' in which the Planck scale, $M_P=10^{19} \GeV$ in
conventional units, is set to $1$. If one knows the dimensions of
a given quantity, it is easy to restore these factors. The fields
$z^i$ conventionally have dimension $[M_P]$ (this is chosen to
make the action $\int |\p z|^2$ dimensionless).  The scalar
potential $V$ and the cosmological constant conventionally have
dimension $[M_P^4]$, while the superpotential $W$ has dimension
$[M_P^3]$.

\subsection{An example from string theory}

Let us describe a simple example of an effective supergravity
theory, which is known to arise from string theory, following the
work of Giddings, Kachru and Polchinski \cite{GKP}.  Further
details can be found in \cite{AD}.

One starts with the IIb superstring theory, and takes the $9+1$
space-time dimensions to be topologically $\R^{3,1} \times X$,
where $\R^{3,1}$ is four-dimensional Minkowski space-time, and $X$
is a Calabi-Yau manifold, a three complex dimensional compact
K\"ahler manifold with zero first Chern class.  It can be shown
that $\dim H^{3,0}(X,\C)=1$ and that the holomorphic three-form
$\Omega$ is nowhere vanishing on $X$.  By Yau's theorem, $X$
admits a Ricci flat metric, so this space-time solves Einstein's
equations.

Furthermore, the moduli space of Ricci flat metrics is isomorphic
to the moduli space of complex structures on $X$, times a
complexified K\"ahler cone. After compactification, this moduli
space forms a factor in the supergravity configuration space $M$,
and each point in $M$ is a possible compactification.  The
K\"ahler metric is simply the Weil-Peterson metric on the moduli
space (the natural metric on the space of metrics).

There is a natural line bundle $\CL$ associated to the K\"ahler
metric. As a bundle over complex structure moduli space, it has a
simple geometric description: it is the Hodge line bundle
$H^{3,0}(X,\C)\rightarrow M$ in which the holomorphic three-form
$\Omega$ takes values. For more about the geometry associated to
this situation, see \cite{Strom,FreedSpec}.

IIb superstring theory contains one more complex scalar field, the
so-called ``dilaton-axion''.  It parameterizes another factor in
$M$, which is the upper half plane with the constant negative
curvature metric. Though approximate, it is standard to take the
metric on $M$ to be a direct product of this metric, with the
Weil-Peterson metric.

A simple example of a section of $\CL$ is a period of $\Omega$.
The superpotentials are the following linear combinations of
periods:
\begin{equation}\label{eq:gvw}
W = \int_X \Omega \wedge \left( F^{(1)} + \tau F^{(2)} \right) ,
\end{equation}
where $\tau$ is the dilaton-axion and $F^{(1)}$ and $F^{(2)}$ are
independently chosen elements of $H^3(M,\Z)$.

This superpotential describes the contribution to the effective
potential due to a ``gauge field strength'' or ``flux'' $F$.  As a
simple indication of this, we note that the formula
(\ref{eq:defV}) implies that $V$ is quadratic in $F$, as is true
for the energy of a magnetic field in Maxwell's theory, and as is
true in supergravity. The standard argument for this
superpotential \cite{GKP} proceeds as follows.  First, one can
show that the critical points $DW=0$ are points in moduli space at
which the form $F^{(1)} + \tau F^{(2)}$ is purely in
$H^{2,1}(X,\C)\oplus H^{0,3}(X,\C)$.  Second, this condition can
be shown to imply that we are at a supersymmetric vacuum.
Finally, (\ref{eq:gvw}) is the unique superpotential with these
properties.

Thus, in this class of compactifications, we obtain a family of
superpotentials, each of which is a linear combination of a finite
basis of sections, the periods of the holomorphic three-form.  To
count vacua in this family of theories, we must count all the
critical points of all of these sections which are allowed
physically.  This could be done by finding the expected number of
critical points of a random section from this class, taken from an
appropriate distribution, and multiplying by the number of
distinct sections.  Thus we have reformulated the physical problem
as a problem in the statistics of holomorphic sections.

In fact, there is a physically well motivated choice for the
ensemble of these sections, which we discuss in detail in
\cite{AD} and will return to in future work.  It consists of the
superpotentials (\ref{eq:gvw}) satisfying the constraint
\begin{equation}\label{eq:length}
\int_X F^{(1)} \wedge F^{(2)} = L
\end{equation}
for some $L\in \Z$. Each superpotential is taken with weight $1$,
to obtain the total number of critical points.

The condition (\ref{eq:length}) sets the overall scale of $W$, and
is analogous to the ``spherical ensemble'' of \cite{BSZ2}, of
sections with a coefficient vector of unit length.  While the
coefficients here must be integers, in the limit
$L\rightarrow\infty$, one expects the sections to be uniformly
distributed on the constraint surface, and thus it should be a
good approximation to neglect the quantization condition; the
spherical or Gaussian ensembles should provide the large $L$
asymptotics for the numbers of physical vacua.

There are further subtleties in making all of this precise, which
we hope to return to in future work.  For one thing, the quadratic
form in (\ref{eq:length}) is indefinite.  The reason this still
defines a ``spherical ensemble'' is a slightly subtle argument
which shows that upon restriction to the subspace of sections with
a critical point at a chosen point in $M$, the form is positive
definite.

In any case, this discussion should convince the reader that the
problem of finding critical points of Gaussian random sections, is
remarkably close to actual problems arising in string theory.

\subsection{Further physical questions about vacua}

Suppose we could find the vacua of the theories we just described,
or of other compactifications of string/M theory: what physical
questions would we like to answer?  Let us discuss questions which
can be answered with the data $(M,K,W)$.

The most basic question is to count the supersymmetric and
nonsupersymmetric vacua, or just the metastable ones.  The
simplest number to obtain is the ``supergravity index,'' which
counts critical points with a weight $\pm 1$ as follows: if $W=0$
at the critical point, the weight is $+1$, while if $W\ne 0$, the
weight is the Morse index. This number is topological for $M$
compact and nonsingular, and more generally is the integral of a
topological density.  Thus, it would be useful to obtain estimates
for the other vacuum counts, in terms of this index, possibly
under conditions such as bounds on curvature and its derivatives.

Besides counting the vacua, we might try to get a picture of their
distribution in the configuration space $M$, by defining a measure
whose integral over a region $R\subset M$ counts vacua within that
region.  To be precise, denote the candidate supergravity theories
as $T_a$, and within each of these, denote the critical points as
$z_i$; the vacuum distribution is then
$$
d\gamma[z_i] = \sum_{T_a} \sum_i \delta_{z_i}
$$

Other distributions over vacua can be defined similarly.  Let
$A_a$ be a function on $M$ in a given theory $T_a$, for example
the cosmological constant or supersymmetry breaking scale.  We
then define its distribution as
$$
d\gamma[A] = \sum_{T_a} \sum_i \delta_{A_a(z_i)} .
$$

A basic question about the supersymmetric vacua, is their
distribution of cosmological constants $\Lambda=-3e^K|W|^2$, and
especially the distribution near zero.  It would be particularly
interesting to find the distribution for flux superpotentials with
integer coefficients.

For nonsupersymmetric vacua, one would like the joint distribution
of cosmological constant $\Lambda$ and supersymmetry breaking
scale $M_{susy}$, ideally just for the metastable vacua.

Finally, one would like to consider more complete definitions of
stability. In particular, a vacuum with $\Lambda=\Lambda_1>0$ can
tunnel or decay to another vacuum with cosmological constant
$\Lambda_2$ satisfying $\Lambda_1>\Lambda_2\ge 0$, at a rate
roughly given by
$$\exp -\int dz \sqrt{V(z)-g_{i\bar j}(z) \dot z^i \dot {\bar z}^{\bar
j}}   ,
$$
{\it i.e.} the exponential of an action, integrated along an
action-minimizing trajectory between the two vacua. This formula
is somewhat simplified, and more precise treatments can be found
in \cite{Cole,Banks,KKLT}, but serves to illustrate the problem.

The total decay rate for a vacuum, is then the sum of this rate
(and, possibly the rate for other decay processes), over all
candidate target vacua. This consideration leads to a constraint
which the vacuum describing our universe must satisfy: its decay
rate should be smaller (hopefully, far smaller) than the inverse
of the known time since the Big Bang, about $10^{10}$ years.
Translated into Planck units, this is about $10^{-60}$. Now in
cases studied so far \cite{KKLT}, the decay rate to any single
target vacuum is far smaller than this, around $10^{-100}$, but it
is conceivable that for $M$ of high dimension, summing the rate
over a large number of targets would lead to an interesting
constraint.

\section{\label{CPHS}Critical points of  holomorphic sections}

We begin the mathematical discussion with the definition of
critical points of a holomorphic section $s \in H^0(M, L)$
relative to  a connection $\nabla $ on $L$.
 We recall that a smooth  connection is a linear map
$$\nabla: \ccal^\infty(M, L) \to \ccal^\infty(M, L \otimes T^*)$$
satisfying $\nabla f s = df \otimes s + f \nabla s$ for $f \in
\ccal^\infty(M)$.  Choosing a local frame
$e_L$ of the line bundle
$L$, we let
\begin{equation}\label{K} K(z)=-\log|e_L(z)|_h^2\;.\end{equation}
The Chern connection
$\nabla =
\nabla_h$ is given by
\begin{equation}\label{Chern} \nabla (f\,e_L) =
(df -f\d K)\otimes e_L\;,\end{equation} i.e., the connection
1-form (with respect to $e_L$) is $-\d K$. We denote the curvature
of $h$ by \begin{equation}\label{Theta} \Theta_h=-d\d K = \ddbar
K\;.\end{equation} (Thus, a positive line bundle $(L,h)$ induces
the \kahler form $\om=\frac i2 \Theta_h=\frac i2 \ddbar K$ with
{\it \kahler potential\/} $K$.)
By \eqref{Chern},
 $\nabla'' s = 0$ for any holomorphic section $s$ where $\nabla =
\nabla' + \nabla''$ is the splitting of the connection into its $L
\otimes T^{*1 ,0}$, resp.\ $L \otimes T^{*0, 1}$ parts.

\begin{defin} \label{DEFCRIT}  Let $(L,h) \to M$ be a holomorphic
line bundle over a complex
manifold, equipped with its Chern connection $\nabla=\nabla_h$. A
critical point of a holomorphic section $s \in H^0(M, L)$ with
respect to $\nabla$  is defined to be a point $z\in M$ where
$\nabla s (z) = 0$, or equivalently $\nabla' s (z) = 0$.   We
denote the set of critical points of $s$  by $Crit^{\nabla}(s)$.
\end{defin}
\medskip

It is important to understand that  the set of critical points
$Crit^{\nabla}(s)$ of $s$, and even its number $\#
Crit^{\nabla}(s)$, depends on $\nabla = \nabla_h$  (or
equivalently on the metric $h$). According to \eqref{Chern},
 the critical point
condition in the local frame, $s =
f e_L$, reads:
\begin{equation} \partial f = f \d K \iff   \partial
\log f =  \partial K\;. \end{equation} As mentioned in the
introduction, this  is a real $\ccal^\infty$ equation, not a
holomorphic one since  $\nabla s \in \ccal^\infty(M, L \otimes
T^{*1, 0})$ is a smooth but not  holomorphic section and
consequently does not
always have positive intersection numbers with the zero section.
Heuristically, the number of critical points reflects the degrees
of both $f$ and of $K$ and the expected number of critical points
should be large if the `degree' of $K$ is large.

An essentially  equivalent  definition in the case of a Chern
connection  is to define a critical point as a point $w$ where
\begin{equation}\label{critfn} d |s(w)|^2_h = 0. \end{equation}
Since
$$d |s(w)|^2_h  = 0 \iff\\[10pt]0=\d |s(w)|^2_h  =
h_w (\nabla' s(w), s(w) ) $$ it follows that \eqref{critfn} is
equivalent to $\nabla' s(w)=0$ as long as $s(w) \not= 0$. So the
critical point condition \eqref{critfn} gives the union of the
zeros and critical points of the section $s$. Another essentially
equivalent critical point equation which puts the zero set of $s$
at $-\infty$ is
\begin{equation} d
\log |s(w)|^2_h = 0.
\end{equation} This is the equation studied by Bott \cite{B}
in his Morse-theoretic proof of the Lefschetz hyperplane theorem,
which is based on the observation that the Morse index of any such
critical point is at least $m$.  We shall use this observation to
study the Morse index density in \S \ref{s-morse}, where we note that the
critical point theory of holomorphic sections at non-singular critical
points is truly just the real Morse theory of the function $\log
|s(z)|_h^2$.

We also note that the classical notion (cf. \cite{AGV, M} ) of
critical point of a holomorphic function $f(z_1, \dots, z_m)$ on
$\C^m$, i.e. a point $w$ where
\begin{equation} \label{CLASSICAL} \frac{\partial f}{\partial z_1} (w) = \cdots
\frac{\partial f}{\partial z_m} (w) = 0 \end{equation}  can be
viewed as a connection  critical point equation in the sense of
Definition \ref{DEFCRIT} but with a   {\it meromorphic
connection\/} rather than  smooth Chern connection. That is, the
derivatives $\frac{\partial f}{\partial z_j}$ on $\C^m$ define a
meromorphic connection
  on the line bundles
$\ocal(N) \to \CP^m$ with poles at infinity. Unlike the case of
smooth connections, the  critical point theory with respect to
meromorphic connections is entirely a holomorphic theory.  The
critical points of a generic section in the sense of
(\ref{CLASSICAL})  all have topological index $+1$, and hence the
number of critical points is a topological quantity depending on
the polar variety of the meromorphic connection and the Chern
classes of $M$ and $L$. This is in contrast to the case of a
smooth connection, where the critical points of a generic section
may have topological index $-1$ as well as $+1$, and their number
depends on the section. As mentioned in the introduction (in the
case of curves), the average number of critical points in the
sense of Definition \ref{DEFCRIT} is greater than the almost sure
number in the classical sense.

The theory of critical points of holomorphic functions (cf.
\cite{AGV, M} ) is concerned with the singularities of the
hypersurface $f(z) = f(z_0)$ at a critical point $z_0$. The
function $g(z) = f(z) - f(z_0)$ has a singular point at $z_0$,
i.e.  $g(z_0) = \nabla g(z_0) = 0$. The same notion of singular
point applies to Definition \ref{DEFCRIT} for  holomorphic
sections. We note that generic holomorphic sections and generic
polynomials have no singular points. Those which do form the
discriminant locus $\dcal \subset H^0(M, L)$. In physics
terminology, singular points are known as Minkowski vacua.
the statistics of singular points are quite different from those of
critical points, and in particular $\dcal$ is a nonlinear
subvariety of $H^0(M, L)$ and  does not carry Gaussian measures.

\subsection{Hessians at a critical point }\label{HESS1}

There are three versions of the Hessian of $s$ at a critical point
which play a role in this paper. In this section, we define them
and explain the relations between them.

The first version of the Hessian of $s$ is
\begin{equation}\label{firstH} D \nabla s (z_0)  \in (T^{*2, 0}
\oplus T^{*1,1})
\otimes L,\;\;\; (\nabla s(z_0) = 0)\;,\end{equation} Here, $D$ is
an auxiliary connection on $T^*M \otimes L$. As is well known, $D
\nabla s(z_0)$ at a critical point is independent of the choice of
$D$. This Hessian will be part  of the jet map
defined in (\ref{JCAL}).

The second version is the `vertical part' $D^v \nabla s$ of the
derivative of the section $\nabla s: M \to T^{* 1,0} \otimes L$
with respect to a connection $D$ on $T^{*1,0} \otimes L$. For
lack of a standard term, we refer to it as the {\it complex
Hessian\/} of $s$. This complex Hessian is the Hessian whose
determinant appears in the statement of Theorem \ref{KNcrit1}. It
is defined as follows:  From an invariant point of view, the
connection gradient $\nabla s$ defines a  section
\begin{equation} \label{NABLAMAP} \nabla s: M \to T^{* 1, 0}
\otimes L. \end{equation}
 We define $D^v \nabla s$ to be  the vertical part of the derivative of
(\ref{NABLAMAP}) with respect to $D$. At a critical point $D^v
\nabla s (z)$ is independent of the choice of the connection $D$.
(The full derivative of $\nabla
s$ maps $TM$ to $T(T^{* 1, 0}
\otimes L)$, which has real dimension $4m$, while $D^v \nabla s$ maps
$TM$ to the vertical tangent space $T^v(T^{* 1, 0}
\otimes L)\approx T^{* 1, 0}
\otimes L$.)

 To define and compute  the various Hessians, we introduce  local
coordinates and an adapted frame in the following sense:

\begin{defin}\label{ADAPTED} Let $\nabla$ be the Chern connection
on a Hermitian holomorphic line bundle $(L,h)\to M$. Let $e_L$ be
a local frame (non-vanishing holomorphic section) of $L$ in a
neighborhood of $z_0\in M$, and let $K$ be the local  curvature
potential given by \eqref{K}. We say that $e_L$  is adapted to
$\nabla$ to order $k$ at $z_0$ if all pure holomorphic
derivatives of $K$ of order $\leq k$ vanish at $z_0$ (and thus
the pure anti-holomorphic derivatives also vanish). In particular,
the connection form vanishes at $z_0$.
\end{defin}

We then write \begin{equation}\label{vj} \nabla s = \sum v_j\,
dz_j\otimes e_L
\;,
\qquad v_j=\frac{\d f}{\d z_j} -f
\frac{\d K}{\d z_j}\;.\end{equation} We fix a point $z_0\in M$ and
choose an adapted local frame (of order 2) at $z_0$ as well as
local normal holomorphic coordinates $z_1, \dots, z_m$ at $z_0$
(i.e., the connection form on $T_M$ also vanishes at $z_0$ in these
coordinates).

We then define  linear functionals $ H'_{jq}, H''_{jq}$
(depending on our choice of coordinates and frame) on the space
$H^0(M,L)$ by:
\begin{equation}\label{Hjq}
D'\nabla' s(z_0)=\sum_{j,q}H'_{jq} dz_q\otimes
dz_j\otimes e_L,\qquad D''\nabla' s(z_0)=\sum_{j,q}H''_{jq}
d\bar z_q\otimes dz_j\otimes e_L\,.\end{equation}

To obtain formulas for the matrices  $H'=\big( H'_{jq}\big)$,
$H''=\big(H''_{jq}\big)$, we recall from
\eqref{K} that
\begin{equation}\label{Kpot}|e_L(z)|_h^2 =
e^{-K(z)}\;,\end{equation} and thus for a section $s=fe_L\in
H^0(M,L)$, we have by \eqref{Chern}:
\begin{equation}\label{covar}
\nabla s=\sum_{j=1}^m\left( \frac{\d f}{\d z_j} -f \frac{\d K}{\d
z_j}\right)dz_j \otimes e_L=\sum_{j=1}^m e^{K}\frac{\d }{\d
z_j}\left(e^{-K}\, f\right)dz_j \otimes e_L\;.
\end{equation}
Differentiating \eqref{covar}, we then obtain:
\begin{eqnarray}H'_{jq} &=& \frac{\d^2 f}{\d
z_j\d z_q}{(z_0)}\;,\label{H'}\\
    H''_{jq} &=& -\left.f \frac{\d^2 K}{\d z_j\d\bar
z_q}\right|_{z_0}=-f(z_0)\Theta_{jq}\,,\quad
\Theta_h(z_0)=\sum_{j,q}\Theta_{jq}dz_j\wedge d\bar
z_q\;.\label{H''}\end{eqnarray}

Thus, the standard Hessian $D\nabla s$ (see \eqref{firstH}) is
given in our adapted coordinates and normal frame by the
$m\times 2m$ matrix $\big(\,H'\ H''\,\big)$, where $H'
$ is a (complex-valued) symmetric matrix, and $H''=-f(z_0)\Theta$,
where
$\Theta$ is the curvature matrix
$\big(\Theta_{jq}\big)$.

To describe the complex Hessian $D^v\nabla s$, we begin by writing
$z_q=x_q+iy_q$ and
$v_j=\sigma_j+i\tau_j$ so that the real Jacobian matrix  (at
$z_0$) of
$\nabla s$ with respect to the variables $\sigma_j,\tau_j$ and
$x_q,y_q$ and the local frame $e_L$ is
\begin{equation}\label{realH}\begin{pmatrix}\ \left(\frac
{\d \sigma_j }{\d x_q} \right) &\left(\frac {\d \sigma_j}{ \d y_q}
\right)\\[10pt]
\left(\frac {\d \tau_j}{\d x_q}\right) &\left(\frac {\d\tau_j}{\d
y_q}
\right) \
\end{pmatrix}\;.\end{equation}
But if we instead compute the Jacobian of $\nabla s$ with respect
to the variables $v_j, \bar v_j$ and $z_q,\bar z_q$, we obtain the
 matrix
\begin{equation}\label{Hmatrix}H^c:=\begin{pmatrix} \left(\frac
{\d v_j }{\d z_q} \right) &  \left(\frac {\d v_j}{ \d
\bar z_q}
\right)\\[10pt]
\left(\frac {\d \bar v_j}{\d
z_q}\right) & \left(\frac {\d \bar v_j}{\d\bar z_q}  \right)
\end{pmatrix}
=\begin{pmatrix} H' &H''\\[6pt] \overline{H''} &\overline{H'}
\end{pmatrix} =\begin{pmatrix} H' &-f(z_0)\Theta\\[8pt]
-\overline{f(z_0)\Theta} &\overline{H'}
\end{pmatrix}\;.
\end{equation} Thus the complex Hessian is represented by the
matrix
$H^c$.

In invariant terms,
 at a critical point $\nabla
s(z_0)=0$, we may express $D^v \nabla s(z_0)$ as the matrix
\begin{equation} \label{Xmatrix} D^v \nabla s (z_0) =
 \left( \begin{array}{cc}
\Hess_{hol}(\frac{\partial}{\partial z_j},
\frac{\partial}{\partial z_k}) s(z) &
{\Theta(\frac{\partial}{\partial
z_j},\frac{\partial}{\partial \bar{z}_k} )
s(z) }  \\ &  \\
\overline{\Theta(\frac{\partial}{\partial
z_j},\frac{\partial}{\partial
\bar{z}_k} ) s(z)} & \overline{\Hess_{hol}(\frac{\partial}{\partial
z_j}, \frac{\partial}{\partial z_k}) s(z))}
\end{array} \right).
\end{equation}
relative to a basis of coordinate vector fields of local
holomorphic coordinates.  Here, the `holomorphic Hessian'
$\Hess_{hol}$ of $s$ at a critical point is given by
\begin{equation}\label{holH}\Hess _{hol}(Z,W) s =
\nabla_{\wt Z}\nabla_{\wt W} s(z_0)\,\qquad Z,W\in
T^{1,0}_{z_0}\,,\end{equation} where $\wt Z,\wt W$ are local
holomorphic vector fields taking the values $Z,W$, respectively,
at $z_0$. Indeed, \eqref{holH} is clearly independent of the
choice of $\wt Z$. Since the curvature $\Theta$ is of type
$(1,1)$,
$$ (\nabla_{\wt Z}\nabla_{\wt W}s-
\nabla_{\wt W}\nabla_{\wt Z}s)(z_0)=  \big(\nabla_{\wt
Z}\nabla_{\wt W}- \nabla_{\wt W}\nabla_{\wt Z}- \nabla_{[\wt Z,\wt
W]}\big)s\,(z_0) = \Theta(Z,W)s(z_0)=0\,,$$ it follows that $\Hess
_{hol}(Z,W) =\Hess _{hol}(W,Z)$, which is also independent of the
choice of $\wt W$.

The off-diagonal terms are the mixed' Hessian,  given by
$$\Hess _{mixed}(Z, \overline W)(s) = \nabla_{\overline
W}\nabla_{ Z}s(z_0).$$ Since
$$\Theta(Z,\overline W)s(z_0) =\big(\nabla_{ Z}\nabla_{\overline
W}- \nabla_{\overline W}\nabla_{ Z}- \nabla_{[Z,\overline
W]}\big)s\,(z_0) = -\nabla_{\overline W}\nabla_{ Z}s(z_0)$$ (here
we dropped the $\ \wt{\ }\,$), the mixed Hessian equals
$$\Hess _{mixed}(Z, \overline W)(s)  =-\Theta(Z,\overline W)s(z_0).$$

The third version is the usual Hessian of $\log |s|^2_h$
at a critical point. This version  will be important in our
discussion of Morse indices in \S \ref{s-morse}. With respect to the
basis $\{dz_j,d\bar z_j\}$, it is given at a critical point $z_0$ by
\begin{eqnarray}  &&\begin{pmatrix}
\left(\frac{\d^2}{\d z_j \d z_q}\log |s|^2_h\right) &
\left(\frac{\d^2}{\d z_j \d\bar z_q}\log |s|^2_h\right) \\[10pt]
\left(\frac{\d^2}{\d\bar z_j \d z_q}\log |s|^2_h\right) &
\left(\frac{\d^2}{\d\bar z_j \d\bar z_q}\log |s|^2_h\right)
\end{pmatrix}\nonumber \ =\ \begin{pmatrix}
\left(\frac 1f \frac{\d^2 f}{\d z_j \d z_q}\right)  &
\left(-\frac{\d^2 K}{\d z_j \d\bar z_q}\right) \\[10pt]
\left(-\frac{\d^2 K}{\d\bar z_j \d z_q}\right) &
\left(\frac 1{\bar f} \frac{\d^2 \bar f}{\d\bar z_j \d\bar z_q}\right)
\end{pmatrix}\\[5pt] &&\hspace{1in} =\  \begin{pmatrix} \frac 1
{f(z_0)}H' &-\Theta\\[8pt] -\overline{\Theta} &{\frac 1
{\bar f(z_0)}\overline{H'}}
\end{pmatrix} \ = \ \begin{pmatrix} f(z_0)\inv &0\\ 0 &\bar
f(z_0)\inv\end{pmatrix} H^c\;.
\label {H3}\end{eqnarray}
Note that the matrix \eqref{H3} is not Hermitian.  In \S
\ref{s-morse}, we use a Hermitian version of \eqref{H3}
obtained by conjugating  the real Hessian of $\log |s|^2_h$ by a unitary
matrix; the resulting Hermitian matrix \eqref{matrix}
contains the entries of
\eqref{H3}, re-arranged and with constant factors.

\section{A density formula for zeros}\label{s-proof}

We now begin the  study of the statistics of critical points of
random sections $s \in \scal \subset  H^0(M, L)$ with respect to a
complex Gaussian measure $\gamma$. We recall that a complex
Gaussian measure is induced by a choice of Hermitian inner product
$\langle, \rangle$ on $\scal \subset H^0(M, L)$:

\begin{defin} \label{GM}  We  define the Gaussian measure associated to $(\scal,
\langle, \rangle)$ by
\begin{equation}\label{gaussinner} d\gamma(s)=\frac 1
{\pi^d}e^{-\|c\|^2} dc\;,\qquad  s=\sum _{j=1}^d
c_je_j,\end{equation} where $dc$ is Lebesgue measure and $\{e_j\}$
is an orthonormal basis for $\scal$ relative to $\langle,
\rangle$.
\end{defin}

We  denote the expected
value of a random variable $X$ on $\scal$ with respect to $\gamma$
by ${\bf E}_\gamma X$ or simply by ${\bf E}X$ when $(\scal, \ga)$ are
understood. We recall from \eqref{CSH}--\eqref{DENSITY} that the expected
distribution of critical points of
$s
\in
\scal \subset H^0(M, L)$ with respect to  $(\scal, \gamma,\nabla
)$ is the measure $\K_{\scal,\gamma,\nabla}=\E_{\ga} C^\nabla_s$ on
$M$, where
\begin{equation*}   C^\nabla_s =\sum_{z\in Crit^{\nabla}(s)}
\delta_{z}\;,
\end{equation*} where $\delta_{z}$ is the Dirac point mass at
$z$. Thus,\begin{equation}
\left(\K_{\scal,\gamma,\nabla},\phi\right) =
\int_{\scal} \left[\sum_{z: \nabla s(z) = 0} \phi(z) \right]\,
d\gamma(s). \end{equation}

\subsection{Covariance kernel}

The crucial invariant of a Gaussian measure is its covariance or
two-point kernel:

\begin{defin} \label{COVKER} The two-point kernel of a Gaussian measure $\gamma$ defined
by $(\scal, \langle, \rangle)$ is defined  by
$$\Pi_{\scal}(z, w) = {\bf E}_{\scal}(s(z) \otimes \overline{s}(w)) \in
L_z \otimes \overline{L_w}.  $$ \end{defin}

Here $\overline L$ denotes the complex conjugate of the line bundle $L$
(characterized by the existence of a conjugate linear bijection
$L \buildrel{\approx}\over {\to}\overline L,\ v\mapsto \bar
v$).  As is well-known and easy to see,
$\Pi_{\scal}$ can be written in the form
$$\Pi_\scal(z,w)=\sum_{j=1}^n s_j(z)\otimes \overline{s_j(w)}\;,$$
where $\{s_1,\dots s_n\}$ is an orthonormal basis for $\scal$ with
respect to the inner product $\langle, \rangle$ associated to the
Gaussian measure $\gamma$. Indeed,
\begin{equation}\label{expect} \E\left(s(z)\otimes \overline
{s(w)}\right) = \E\left(\sum_{j,k=1}^nc_j\overline{c_k}\, s_j(z)
\otimes \overline {s_k(w)}\right) =  \sum_{j=1}^n s_j(z)\otimes
\overline{s_j(w)}\;,\end{equation} since the $c_j$ are independent
complex (Gaussian) random variables of variance 1.

In the  case of a Hermitian line bundle,  the   two point kernel
of the Hermitian Gaussian measure  is the \szego kernel of $(L,
h)$, i.e.   the orthogonal projection $ \Pi_{\scal, h, V}:
\lcal^2(M, L) \to \scal$ with respect to the inner product
(\ref{HIP}).

\subsection{Expected density of random discrete zeros}

 The
expected density of critical points may be regarded as the
expected density of zeros of random sections in the subspace
$\nabla H^0(M, L) \subset  \ccal^\infty(M, T^{*1,0} \otimes L)$. In
this section, we prove a general formula (Theorem \ref{density})
for the density of zeros of random sections which applies to this
subspace and which  will be used to prove Theorem \ref{KNcrit1}.
It may be derived from the rather general and abstract Theorem~4.2
of \cite{BSZ3}. However, that theorem gives the $n$-point
correlation of zeros of several random sections  in all
codimensions, while here we consider only the density (or
``1-point correlation") in the full codimension case where the
zeros are discrete.  This is both simpler than the general setting
in \cite{BSZ3} and also involves some special features not quite
covered there. To  make the paper more self-contained, we give a
derivation from scratch of the density formula for discrete zeros
that arises from \cite[Theorem~4.2]{BSZ3}. In \S\ref{s-alt}, we
give an alternate approach to the proof which is closer to
\cite{AD}.

The general set-up in \cite{BSZ3} involves $1$-jets of sections of
a real vector bundle $V$ over a smooth manifold $M$.  (We shall
later apply our formula to the case where $V=T^{*1,0} \otimes L$
is complex, but the sections $\nabla s\in\scal$ are not
holomorphic.) For simplicity of exposition, we will endow $V$ with
a connection $\nabla$ and an inner product $h$, and we will
endow $TM$ with a Riemannian metric and a volume form $d\vol_M$.
The result of Theorem \ref{density} below is independent of these
choices of connection and metric.

 Let $\scal \subset \ccal^\infty (M, V)$ be a
finite-dimensional subspace of smooth sections  and consider the
{\it jet maps}
$$J^1_z:  \scal \to J^1(M,V)_z\;,
\qquad z\in M\;, $$
where
$J^1(M, V)$ denotes the vector bundle of
$1$-jets of sections of $V$, and $J^1_z(s)$ is the 1-jet at $z\in M$
of a section $s\in\scal$.  Recall that we have the canonical vector bundle
exact sequence
\begin{equation} \label{exactseq} 0\to T_M^*\otimes V \to  J^1(M,
V)\buildrel{\epsilon}\over{\to} V \to 0\;,\end{equation}
where $\epsilon$ is the evaluation map.

The connection $\nabla$ on $V$ gives a splitting of
\eqref{exactseq},
\begin{equation}\label{splitting}(\ep,\nabla):J^1(M,V)
\buildrel\approx\over\to  V\oplus (T_M^*\otimes V)\,,\qquad J^1_z(s)\mapsto
(s(z),\nabla s(z))\,.\end{equation} We shall identify
$V\oplus (T_M^*\otimes V)$ with the space $J^1(M,V)$ of 1-jets via
\eqref{splitting}. Given a Gaussian measure $\gamma$ on $\scal$ and a point
$z\in M$, we consider the pushforward measure
\begin{equation}\label{JPD}\D_z:=(J^1_z)_*
\gamma\;,\end{equation} which is called the  {\it joint
probability distribution\/} of
$\gamma$.  Since the jet map $J^1_z$ is Gaussian, the joint
probability distribution $\D_z$ is
likewise Gaussian.

In the application we have in mind, $V = T^{*1,0} \otimes L$,
$\scal = \nabla H^0(M, L)$ and
\begin{equation} \label{JCAL} J^1_{z}:\nabla H^0(M, L) \to J^1(M,T^{*1,0}
\otimes L) \approx
\left(T^{*1,0} \oplus [T^{*1,0}\otimes  T^{*1,0}] \oplus  [T^{*0,1}\otimes
T^{*1,0}]\right)_z
\otimes L_z.
\end{equation} A complication  arises (when
$\dim M>1$) in  that the range of $J^1_z$ is a proper subspace of $J^1(M,
V)_z$. Indeed, in terms of normal coordinates, $J^1(M,T^{*1,0}
\otimes L)$ can be identified with the space of triples $(v,H',
H'')$, where $v\in\C^m$ and $H',H''$ are complex $m\times m$ matrices,
while the range  of
$J^1_z$ consists only of those triples
where $H'$ is a complex symmetric matrix and $H''=x\,\Theta$,
$x\in\C$ (see (\ref{H'})--(\ref{H''})). Then $\D_z$ becomes a singular
Gaussian measure on $J^1(M, V)$. The results of
\cite{BSZ2} include singular measures, but it is simpler to apply
the results in a way which is better adapted to the subspace
situation.

Hence, returning to our general setup, we assume that the jet map has
the following spanning property:

\begin{defin}\label{spanning}Let  $\scal$ be a linear space of sections
of a $\ccal^\infty$ vector bundle $V\to M$ and let $\jcal^1:M\times\scal \to
J^1(M,V)$ be given by $\jcal^1(z,s)=J^1_z(s)$.  We say that
$\scal$ has  the spanning property with respect to a
sub-bundle
$W\subset T^*_M\otimes V$ if\/ {\rm Image}$\,\jcal^1$ is a sub-bundle of
$J^1(M,V)$ and $\epsilon:{\rm Image}\,\jcal^1 \to V$ is surjective with
kernel  $W$; i.e.,
\begin{equation} \label{SP}0\to W\to \mbox{\rm Image} \,\jcal^1
\buildrel\epsilon\over\to V\to 0 \end{equation} is an exact sequence of
vector bundles.
\end{defin}

The pushforward measure
$\D_z$  of (\ref{JPD}) is then a (nonsingular)  Gaussian measure on
Image$\,J^1_z$. Making the identification
$J^1(M,V)\approx V\oplus(T^*_N\otimes V)$ via \eqref{splitting}, we
have
\begin{equation}\mbox{Image}\,J^1_z \approx V_z\oplus W_z
\;  . \end{equation} We then regard
$\D_z$ as a  Gaussian measure on $V_z\oplus W_z$, and we write
\begin{equation}\label{Dgauss}
\D_z= D(x,\xi;z)\,dx\,d\xi\ \qquad (z\in M,\ x\in V_z,\
\xi\in W_z)\,,
\end{equation} where $dx,\ d\xi$ denote Lebesgue measure on $V_z,\
W_z$ respectively (with respect to our Riemannian metric $G$ on $M$ and
inner product $h$ on $V$).  We
note that
$D(x,\xi;z)$ depends on the choice of metrics, but of course $\D_z$ does
not.

We now assume further that rank$\,V=\dim M=k$, so that by the spanning
property,  the zero sets $Z_s$ of sections
$s
\in
\scal$ are almost surely discrete.
We shall denote by $|Z_s|$ the sum of delta functions at the zeros
of $s$. The following theorem is a special case of Theorem~4.2 in
\cite{BSZ3}:

\begin{theo}\label{density}  Let
$V \to M$ be a $\ccal^\infty$ real vector bundle over a
$\ccal^\infty$ manifold of dimension $k=\mbox{\rm rank}(V)$, and
let
$\scal
\subset \ccal^\infty(M, V)$ be a finite-dimensional subspace with the
spanning property \eqref{SP} with respect to a subspace $W\subset
T^*_M\otimes V$. Let
$\gamma$ be a Gaussian probability measure on $\scal$.   Then
\begin{equation}\label{d1} \E_\ga|Z_s| =\kcal\,d\vol_M \,,\quad
\kcal(z)= \int_{W_z} D(0,\xi;z) \,\|\det\xi\|\, d\xi\,,\end{equation}
where $d\xi$ denotes Lebesgue measure with respect to the metric
on $W_z\subset T^*_{M,z}\otimes V_z$, and where $D(0,\xi;z)$ is
given by \eqref{JPD} and \eqref{Dgauss}. (An explicit formula for
$D(0,\xi;z)$ is given in (\ref{D0}).
\end{theo}

The notation  $\|\det\xi\|$ in \eqref{d1} is defined as follows:
a $V$-valued 1-form $\xi\in (T_M^*\otimes V)_z =
\mbox{Hom}(T_{M,z},V_z)$ induces a $(\det V)$-valued $k$-form
$$\det\xi\in \mbox{Hom}(\det T_{M,z},\det V_z)= \left(\textstyle\bigwedge^k
T_M^*\otimes
\det V\right)_z\,.$$
Then  $\|\det\xi\|$ is the norm on $\det T_{M,z}^*\otimes \det V_z$
induced from the metrics on $M$ and $V$. To describe the norm
explicitly,  we write
$$\xi=\sum_{j=1}^k  \xi_j\otimes e_j, \qquad \xi_j\in T^*_z,$$
where $\{e_1,\dots,e_k\}$ is an orthonormal basis for $V_z$. Then
\begin{equation}\label{det}\|\det\xi\|= \|\xi_1\wedge \cdots\wedge \xi_k\|=
\left|\frac{\xi_1\wedge \cdots\wedge
\xi_k}{d\vol_M}\right|\;.\end{equation}

\begin{rem} We note that $D(0,\xi;z)$ is independent of the choice of the
connection $\nabla$ (see \cite[p.~371]{BSZ2}). (It does depend on
the choice of metric on $V$, but the reader can easily check that
$D(0,\xi;z) \,\|\det\xi\|\, d\xi\, d\vol_M(z)$ defines a measure
on $W$ that is independent of metrics and volume forms.)
\end{rem}

\subsubsection{Zeros of sections of complex vector bundles.}
Now let $V\to M,\ \scal
\subset \ccal^\infty(M, V)$ be as in Theorem \ref{density}, but let
$V$ be a complex vector bundle of rank $k$ over $\C$. We suppose that
$\dim M=2k$ so that we have point zeros. We may apply
Theorem~\ref{density}, regarding
$V\to M$ as a real vector bundle of rank $2k$.

Then \eqref{d1} holds, but we must properly interpret $\|\det\xi\|$.
To do this, we fix
$z\in M$, and we pick an orthonormal basis $\{e_1,\dots,e_k\}$ of
$V_z$ over $\C$. We then regard $V_z$ as a real vector bundle
endowed with the inner product  having orthonormal basis
$$\left\{\frac 1{\sqrt 2}e_1,\frac i{\sqrt 2}e_1,\dots, \frac 1{\sqrt
2}e_k,\frac i{\sqrt 2}e_k\right\}\;.$$  As before, for $\xi\in
(T_M^*\otimes V)_z$, we write
$$\xi=\sum_{j=1}^k \xi_j\otimes e_j=\sum_{j=1}^k\left( \Re\xi_j\otimes e_j
+ \Im\xi_j\otimes ie_j\right), \qquad \xi_j\in T^*_z\otimes
\C\;.$$ Thus we have
\begin{equation}\label{s-cvb}\|\det\xi\|=
2^k\|\Re\xi_1\wedge\Im\xi_1 \wedge\cdots\wedge
\Re\xi_k\wedge\Im\xi_k\| =\|\xi_1\wedge\cdots\wedge\xi_k \wedge
\bar\xi_1\wedge\cdots\wedge\bar\xi_k\|\;.
\end{equation}

\subsubsection {Proof of Theorem \ref{density}.}\label{theproof}
As mentioned above, the theorem is a special case of  Theorem~4.2 in
\cite{BSZ3}.  However, the proof in this case (which is based on the proof in
\cite{BSZ2}) is quite simple, so we present it here.

We can restrict to a
neighborhood $U$ of an arbitrary point $z_0\in M$. Since $\scal$
spans $V$, we can choose $U$ so that there exist
sections
$e_1,\dots,e_k\in\scal$ that form a local frame for $V$ over $U$.
For a section $s\in\scal$, we write $s(z)=\sum_{j=1}^k
s_j(z)e_j(z)$ ($z\in U$) and we let $\tilde s=(s_1,\dots,s_k):U\to
\R^k$.  Since $D(0,\xi;z)$ is independent of the connection, we can
further assume that $\nabla|_U$ is the flat connection $\nabla s =
\sum ds_j \otimes e_j$. Then
\begin{equation}\label{bbb}\|\det\nabla s\| =\sqrt{h} \|d
s_1\wedge\cdots\wedge d s_k\|\,,\end{equation} where
$h=\det\big(h(e_j,e_{j'})\big)$.

We let $\psi_\ep\to\de_0$ be an approximate identity on $\R^k$, and we write
$c=(c_1,\dots,c_k)\in\R^k, \ dc=dc_1\cdots dc_k$. For a
test function $\phi\in\dcal(U)$, we have by \eqref{det} and
\eqref{bbb},
\begin{eqnarray} \int_{\R^k}
\,\psi_\ep(c)\Big(|\tilde s\inv(c)|,\phi\Big)\,dc &=& \int_{\R^k}
\,\psi_\ep(c)\Big[\sum_{\tilde s(z)=c}\phi(z)\Big]\,dc \nonumber\\
&=&  \int_U \,(\psi_\ep\circ \tilde s)\; \phi\, |ds_1\wedge \cdots
\wedge ds_k| \nonumber\\
&=& \int_U \,(\psi_\ep\circ \tilde s)\; \phi\, \|\det\nabla
s\|\,h^{-1/2}\,d\vol_M \,.\label{change}\end{eqnarray}

Integrating (\ref{change}) over $\scal$ and using \eqref{JPD}, we
obtain
\begin{eqnarray} \int_{\R^k}\psi_\ep(c)
\big(\E|\tilde s\inv(c)|,\phi\big)\,dc &=& \int_\scal\int_U
\,(\psi_\ep\circ \tilde s)(z)\, \phi(z)\, \|\det\nabla
s\|_z\,h(z)^{-1/2}\,d\vol_M(z)\,d\ga(s)\nonumber \\
&=& \int_M\int_{W_z}\int_{\R^k}\psi_\ep(c)\,\phi(z)\,
\|\det\xi\|\, D\left(\textstyle\sum
c_je_j,\xi;z\right)dc\ d\xi\,d\vol_M(z),\nonumber\\[-6pt]
\label{E}\end{eqnarray} where the latter equality follows from the fact that
$$(J^1_z)_* (
d\gamma)= D(x,\xi;z)\,dx\,d\xi= D\left(\textstyle\sum
c_je_j,\xi;z\right)\,h(z)^{1/2}\,dc\,d\xi\;.$$

Letting $\ep\to 0$ in \eqref{E}, we obtain
$$\E(|\tilde s\inv(0)|,\phi) = \int_M\int_{W_z}\phi(z)\,
\|\det\xi\|\, D\left(0,\xi;z\right)\,d\xi\,d\vol_M(z)\;.$$
Recalling that $\tilde s\inv(0)=Z_s$, we then obtain \eqref{d1}.\qed

\medskip
\begin{rem} The proof of the analogous result for the case where
rank$\,V<\dim M$ follows the same argument. The only additional ingredient
is Federer's co-area formula, which is used to obtain \eqref{change}; see
\cite{BSZ3}.\end{rem}

\subsection{Description of the joint probability
distribution.}\label{JPDSECT}
We again suppose that $V$ is a complex vector bundle. Recall that the
measure $\D_z$ is the pushforward of the Gaussian measure $\ga$  under
the linear map  $J^1_z$.  Since the
push-forward of a Gaussian measure under a linear map is Gaussian, $\D_z$
is a Gaussian measure on Image$\,J^1_z$.  We now give a formula for $\D_z$
and more importantly, for the conditional Gaussian measure $\D^0_z$ that
appears in our formula
\eqref{d1}.

Let $z\in M$, and choose  orthonormal bases
$\{e_1,\dots,e_k\}$,
$\{ w_1,\dots,w_n\}$ of $V_z,\ W_z$, respectively.  The Gaussian measure
$\D_z$ can be written in the form
\begin{eqnarray*}d\ga_{\De(z)}(v,w)=  \frac 1 {\pi^{k+n} \det
\De(z)}\exp\left[-\left\langle
\De(z)^{-1} \begin{pmatrix}x\\
y\end{pmatrix},\begin{pmatrix}x\\
y\end{pmatrix}\right\rangle\right]\,dx_1\cdots dx_k\,dy_1\cdots dy_n\;,\\
v=\sum_{j=1}^kx_je_j\;,\ w=\sum_{q=1}^ny_qw_q\;.
\end{eqnarray*}
The  covariance matrix $\De(z)$ is given in block form by
\begin{eqnarray} \label{D1}
\Delta(z) =\left(
\begin{array}{cc}
A & B \\
B^{*} & C
\end{array}\right)\,,\quad A=\big[\E(x_j\bar x_{j'})\big],\
 B=\big[\E(x_j\bar y_q)\big],\ C=\big[\E(y_q\bar y_{q'})\big],\\
1\le j,j'\le k,\ 1\le q,q'\le n. \nonumber \end{eqnarray}
Using the
formula for the inverse  of a matrix in block form, we obtain
\begin{equation}\label{D0}
D(0,y;z) =\frac{1}{\pi^{k+n}\det A\det\La}\exp\left(
-{\langle \La^{-1}y,y\rangle}\right)\;,
\end{equation}
where $
\La=C-B^*A^{-1}B$ as in \eqref{Lambda}.

\section{Density formulas: Proof of Theorem
\ref{KNcrit1}}\label{COMPLETION}

We now prove Theorem \ref{KNcrit1} for the ensemble $(\scal,
\nabla, \gamma)$  by applying the zero-density formula
\eqref{d1}--\eqref{det} of Theorem \ref{density}  to the ensemble
$$\scal':=\nabla\scal\subset
\ccal^\infty(M,T^*_M\otimes L)\;,$$  endowed with the Gaussian
probability measure on $\scal'$  induced by $\gamma$.
We assume that $\scal$ has the following property:

\begin{defin} \label{2JET} Let $(L, h) \to M$ be a Hermitian
holomorphic line bundle. We
say that $\scal \subset H^0(M, L)$ has the {\it $2$-jet spanning
property}
 if  the jet maps
$$J_z^2: \scal  \to J^2_{\hol}(M,L)_z$$ are surjective for all $z\in
M$ (where
$J^2_\hol(M,L)$ denotes the vector bundle of 2-jets of holomorphic sections
of
$L$).
\end{defin}

When $L\to M$ is a positive line bundle on a compact complex
manifold $M$, the  surjectivity of $J^2_z$ always holds for
$\scal=H^0(M, L^N)$ when $N$ is sufficiently large, as   a well
known consequence of the Kodaira Vanishing Theorem.

We begin with the following observation:

\begin{lem} Let $(L,h)\to M$ be a Hermitian holomorphic line bundle such
that the Chern connection $\nabla$ had nonvanishing curvature form
$\Theta$.  Suppose that
$\scal\subset H^0(M,L)$ is a linear space of sections with the $2$-jet
spanning property.  Then $\nabla\scal\subset \ccal^\infty(M,
T^{*1,0}_M\otimes L)$ has the spanning property with respect to
$$W:=(S^2T^{*1,0}_M \oplus \C\wt\Theta)\otimes L \subset
(T^{*1,0}_M\otimes T^{*1,0}_M\otimes L)\oplus (T^{*0,1}_M\otimes
T^{*1,0}_M\otimes L)= T^*_M\otimes T^{*1,0}_M\otimes L\;,$$ where
$S^2T^{*1,0}_M \subset T^{*1,0}_M \otimes T^{*1,0}_M$ denotes the
symmetric tensors, and $\wt\Theta$ corresponds to $\Theta$ under
the natural identification $T^{*1,1}_M \approx
T^{*0,1}_M \otimes T^{*1,0}_M$.
\end{lem}

\begin{proof} We begin by describing the relevant random variables
$x,v_j, H'_{jq}, H''_{jq}$ used to describe the  jet map $J^1_{z}$.
Let $z_0\in M$ and choose normal coordinates $\{z_j\}$ and a special
frame
$e_L$ adapted to
$\nabla$  at $z_0$.
Recalling \eqref{vj}--\eqref{Hjq}, we consider the  linear
functionals
$x,v_j, H'_{jq}, H''_{jq}$  on the space
$H^0(M,L)$  given by:
\begin{eqnarray} & {\displaystyle s(z_0)=x\,e_L,\qquad \nabla
s(z_0)=\nabla' s(z_0) = \sum_{j=1}^m v_j dz_j\otimes e_L\,}
\label{X}
\\
& \displaystyle D'\nabla' s(z_0)=\sum_{j,q}H'_{jq} dz_q\otimes
dz_j\otimes e_L,\qquad D''\nabla' s(z_0)=\sum_{j,q}H''_{jq}
d\bar z_q\otimes dz_j\otimes e_L\,.\label{H}\end{eqnarray}

The jet map in local coordinates, using the identification
\eqref{JCAL}, is given by \begin{equation} \label{JONE} J^1_{z_0}=
(v_j, H'_{jq}, H''_{jq})\;. \end{equation} The conclusion is an
immediate consequence of the 2-jet spanning property of $\scal$
and \eqref{H'}--\eqref{H''}.\end{proof}

We recall that the matrices $[H'_{jq}]$ and $[H^{''}_{jq}]$ are
the coordinate representation of  the holomorphic Hessian and
mixed Hessian described in \S \ref{HESS1}, where it was observed
that they form part of the vertical derivative matrix
$H^c$ of $\nabla s$.

\subsection{Density formula and covariance kernel}\label{s-szego}

Following \S \ref{JPDSECT}, we next  compute the joint probability
density using the coordinates $\{H'_{jq}\ (1\le j\le q\le m),\
x\}$ with respect to the basis
$$\{dz_j\otimes  dz_q \otimes e_L|_{z_0}\ (1\le j\le q\le m),\
\Theta_h \otimes e_L|_{z_0}\}$$ of $W_{z_0}$. (Here, in order to
obtain the result as a consequence of Theorem \ref{density} on zero
densities,  we need to assume that the curvature form
$\Theta_h$ does not vanish at
$z_0$.  However, in the general case, the formula follows directly from
the argument in \S
\ref{theproof} using instead the joint probability distribution
$D(v,H',x;z_0)\,dv\,dH'\,dx$.)  The joint probability density
$D(v,H',x;z_0)$ is  Gaussian with covariance matrix  $\Delta(z_0)$
given by:
\begin{eqnarray} \Delta(z_0)&=&\left(
\begin{array}{cc} A & B \\ B^* & C
\end{array}\right)\,, \label{Delta'}\\ A&=& \Big[\E\big
(v_j\overline{ v_{j'}}
\big)\Big]\,,\label{A}\\ B&=& \Big[ \E
\big(v_j\overline{H'_{j'q'}}\big)\quad \E\big( v_j\bar x\big)\Big]
\,,\label{B}\\[8pt] C&=&\left[
\begin{array}{cc}\E \big(H'_{jq}\overline{H'_{j'q'}}\big) &
\E \big(H'_{jq}\bar x\big) \\[8pt]
\E \big(x\overline{H'_{j'q'}}\big) & \E
(|x|^2)\end{array}\right]\,,\label{C}\\ && \qquad\qquad \qquad
1\le j\le m\,, 1\le j\le q\le m\,, 1\le j'\le q'\le m\,.\nonumber
\end{eqnarray}

We now describe how $\De(z_0)$ is given in terms of the covariance
kernel $\Pi_\scal$ of $\scal$ (cf. Definition \ref{COVKER}). It is
in fact simpler to use the local expression for the covariance
kernel in a local  frame (non-vanishing local holomorphic
section). We fix a point $z_0\in M$, and choose a frame
 $e_L$ of $L$ on a neighborhood
$U\subset M$ of $z_0$. We write every section in the form $s
=f\,e_L$.

\begin{defin} \label{xexpect}  The local covariance kernel $F_{\scal}(z, w) \in \ocal (U \times \overline{U})$
 in a frame $e_L$ of $L$
is defined by $$ \Pi_{\scal}(z, z) = F_{\scal}(z, z)\, e_L(z) \otimes
\overline{e_L(z)}\;.$$ Equivalently,
 $$ F_{\scal}(z_0, w_0) =
 \sum_{j} f_j(z_0) \overline{f}_j(w_0)$$
 where $s_j = f_j e_L$ is an orthonormal basis of $(\scal,
 \langle, \rangle).$
\end{defin}

We then have:
\begin{equation}\label{xexpect1}\E(|x|^2)= F_{\scal}(z_0, z_0) =
\sum|f_j(z_0)|^2.
\end{equation}
We emphasize  that both the random variable $x$ of (\ref{X}) and
the formula (\ref{xexpect1}) depend on the choice of frame $e_L$.
It is convenient to introduce an invariant  notation for the local
covariance kernel in the frame $e_L$.  We write (\ref{xexpect1}) as
\begin{equation}\label{xexpect2} \E(|x|^2) = \frac
{\Pi_\scal(z_0,z_0)}{e_L(z_0)\otimes \overline{e_L}(z_0)} =
\rho^\diag_{e_L(z_0)}\,\Pi_\scal\;,
\end{equation}
where $\rho^\diag$ denotes the restriction to the diagonal, and
\begin{equation}\label{rhodiag} \rho_v^\diag G = {G(z_0,z_0)}/(v\otimes \bar v)\in\C\;,\end{equation}
for $ G(z_0,z_0)\in L_{z_0}\otimes \bar L_{z_0},\ v\in L_{z_0}$.

Differentiating \eqref{expect}, we obtain
\begin{equation*}\E\left(\nabla_{z_j}s(z)\otimes \overline
{\nabla_{w_{j'}}s(w)}\right) = \nabla_{z_j} \nabla_{\bar w_{j'}}
\Pi_\scal(z,w)\;,\end{equation*} where we write
$$ \nabla' s=\sum_{j=1}^m dz_j
\otimes \nabla_{z_j} s\,,\quad \nabla'' s= \sum_{j=1}^m d\bar z_j
\otimes \nabla_{\bar z_j} s\;.$$
Hence,
\begin{equation*}\E(v_j\overline{ v_{j'}})
=\rho^\diag_{e_L(z_0)}\nabla_{z_j}
\nabla_{\bar w_{j'}}
\Pi_\scal\;.
\end{equation*} Thus, after repeatedly differentiating
\eqref{xexpect}, the matrices \eqref{A}--\eqref{C} can be
expressed in terms of the covariance kernel and its covariant
derivatives on the diagonal:
\begin{eqnarray}
A&=& \left(\rho^\diag_{e_L(z_0)}\nabla_{z_j}\nabla_{\bar
w_{j'}}\Pi_\scal \right),\label{AN0}\\
B&=&\left[\left( \rho^\diag_{e_L(z_0)} \nabla_{z_j}\nabla_{\bar
w_{q'}}\nabla_{\bar w_{j'} }\Pi_\scal\right) \quad
\left(\rho^\diag_{e_L(z_0)}
\nabla_{z_j}\Pi_\scal\right)
\right] \,,\label{BN0}\\[8pt]
C&=&\left[
\begin{array}{cc}\left(
\rho^\diag_{e_L(z_0)} \nabla_{z_q}\nabla_{z_j}\nabla_{\bar
w_{q'}}\nabla_{\bar w_{j'}}\Pi_\scal \right) &
\left( \rho^\diag_{e_L(z_0)}
\nabla_{z_q}\nabla_{z_j}\Pi_\scal\right)
\\[8pt]
\left(\rho^\diag_{e_L(z_0)} \nabla_{\bar w_{q'}}\nabla_{\bar
w_{j'}}\Pi_\scal\right)& \rho^\diag_{e_L(z_0)}
\Pi_\scal\end{array}\right] \,,\label{CN0}\\
&& \qquad\qquad \qquad 1\le j\le m\,, 1\le j\le q\le m\,, 1\le
j'\le q'\le m\,.\nonumber
\end{eqnarray}
In the above, $A,B,C$ are $m\times m,\,m\times n,\, n\times n$
matrices, respectively, where $n=\half(m^2+m+2)$.

We  pause to obtain  simple local formulas in an adapted frame and
in normal coordinates for $\nabla$. We first replace  each
covariant derivative  by its local expression $\nabla_{z_j} =
\frac{\partial}{\partial z_j} - \frac{\partial K}{\partial z_j} $
in the frame $e_L$ and each $\Pi_{\scal}$ can be replaced by its
local expression $F_{\scal}$. Thus,
\begin{equation} \label{SIMPLER} \E\left(\nabla_{z_j}s(z)\otimes
\overline {\nabla_{w_{j'}}s(w)}\right) = (\frac{\partial}{\partial
z_j} - \frac{\partial K}{\partial z_j}) (\frac{\partial}{\partial
\bar{w}_{j'}} - \frac{\partial K}{\partial \bar{w}_{j'}})
F_{\gamma} (z,w)|_{z = w}.
\end{equation}
Similarly for higher covariant derivatives.

 Thus we have,
\begin{eqnarray}
A&=& \textstyle\left((\frac{\partial}{\partial z_j} - \frac{\partial
K}{\partial z_j}) (\frac{\partial}{\partial \bar{w}_{j'}} -
\frac{\partial K}{\partial
\bar{w}_{j'}}) F_{\scal} (z,w)|_{z = w}\right) ,\label{AN}\\
B&=& \textstyle\left[\left( \frac{\partial}{\partial z_j} - \frac{\partial
K}{\partial z_j})(\frac{\partial}{\partial \bar{w}_{q'} } -
\frac{\partial K}{\partial \bar{w}_{q'}
})(\frac{\partial}{\partial w_{j'} } - \frac{\partial K}{\partial
\bar{w}_{j'}}) F_{\scal}|_{z = w} \right) \quad \left(
(\frac{\partial}{\partial z_j} - \frac{\partial K}{\partial z_j})
F_{\scal}|_{z = w}\right)
\right] \,,\label{BN}\\[8pt]
C&=&\textstyle\left[
\begin{array}{cc} C'& \left( (\frac{\partial}{\partial z_j} - \frac{\partial
K}{\partial z_j})(\frac{\partial}{\partial z_q} - \frac{\partial
K}{\partial z_q}) F_{\scal}|_{z = w} \right)
\\[8pt]
\left( (\frac{\partial}{\partial \bar{w}_{q'}} - \frac{\partial
K}{\partial \bar{w}_{q'}})(\frac{\partial}{\partial \bar{w}_{j'}}
- \frac{\partial
K}{\partial \bar{w}_{j'}}) F_{\scal}|_{z = w}\right) & F_{\scal}(z,z)
 \end{array}\right] \,,\nonumber\\[8pt]
&&\textstyle C'=\left(
(\frac{\partial}{\partial z_q} - \frac{\partial K}{\partial
z_q})(\frac{\partial}{\partial z_j} - \frac{\partial K}{\partial
z_j})(\frac{\partial}{\partial \bar{w}_{q'}} - \frac{\partial
K}{\partial \bar{w}_{q'}})(\frac{\partial}{\partial \bar{w}_{j'}}
- \frac{\partial K}{\partial \bar{w}_{j'}})F_{\scal}|_{z = w}
\right) \label{CN}\\
&& \qquad\qquad \qquad 1\le j\le m\,, 1\le j\le q\le m\,, 1\le
j'\le q'\le m\,.\nonumber
\end{eqnarray}

In the formulas (\ref{AN})--(\ref{CN}), we only take
repeated holomorphic or anti-holomorphic derivatives of the potential $K$.
Hence in an adapted frame of  high order $2$, the matrices simplify to
the ones given in (\ref{AN1})--(\ref{CN1}) when evaluated at $z=w=z_0$.

\subsection{Completion of the proof} From (\ref{D0}), we obtain
\begin{equation}\label{DNm}
D(0,H',x;z_0) =\frac{1}{\pi^{{m+2\choose 2}}\det
A\det\La}\exp\left( -{\langle
\La^{-1}(H'\oplus x),H'\oplus x\rangle}\right)\;,
\end{equation} where $\Lambda$ is given by \eqref{Lambda}.
The formula in Theorem \ref{density}  then yields the expected
density of critical points:
\begin{equation}\label{d2crit}
K(z_0)= \int_{\C^{n}} \|\det H\| D(0,H', x;z_0) \, dH'\, dx\,.
\end{equation}

To complete the proof of Theorem \ref{KNcrit1}, we need the formula
for
$\|\det H\|$.  We now obtain the formula from
\eqref{s-cvb} with $h = 1$ (normal coordinates at $z_0$) and
\begin{eqnarray*}\xi_j  &=&
H'_{j1}dz_1+\dots +H'_{jm}dz_m+ H''_{j1}d\bar
z_1+\dots+H''_{jm}d\bar z_m\,,\\
\bar\xi_j &=& \overline{ H''_{j1}}dz_1+\dots +\overline{
H''_{jm}}dz_m+ \overline{ H'_{j1}}d\bar z_1+\dots+\overline{
H'_{jm}}d\bar z_m\,.
\end{eqnarray*}
By \eqref{s-cvb}, we see that $\|\det H\|$ is the determinant
of the matrix $H^c$ given by (\ref{Hmatrix}):
\begin{eqnarray}\|\det H\| &=& |\det H^c|\ =\
 \left|\det \left[\begin{array}{cc}
H' & -x\Theta\\ -\bar x \bar\Theta &
H'{}^*\end{array}\right]\right|\;. \label{bbbH}\end{eqnarray}
Theorem  \ref{KNcrit1} now follows from \eqref{DNm}--\eqref{bbbH}.
\qed

\begin{rem} Our formulas for the $A,B,C$ matrices differ slightly from
those in \cite{BSZ3} (and in our forthcoming paper \cite{DSZ}).
There, we lift the computations on a positive Hermitian  line
bundle to the associated circle bundle, which amounts to replacing
$\rho^\diag_{e_L(z_0)}$ by $\rho^\diag_{u}$, where
$u=\|e_L(z_0)\|_h\inv e_L(z_0)$.
 The resulting formula for the density is the
same in both approaches, since it is invariant when $\De(z_0)$ is
multiplied by a scalar factor.\end{rem}

\subsection{Proof of Corollary  \ref{KNcrit2}.}

We apply Theorem \ref{KNcrit1} with $\Theta=I$, which is the local
formula for $\Theta$ in normal coordinates;  equivalently
$H''_{jq}=-\de_j^q x$. Then
\begin{eqnarray}\|\det H\| &=& \left|\det \left[\begin{array}{cc}
H' & H''\\ \overline{ H''} & \overline{
H'}\end{array}\right]\right|\ =\
 \left|\det \left[\begin{array}{cc}
H' & -xI\\ -\bar x I & H'{}^*\end{array}\right]\right|\nonumber\\
&=& \left| (\det H')\det (H'{}^* -|x|^2 H'{}\inv)\right|\nonumber\\
&=& \left|\det(H'H'{}^*-|x|^2I)\right|\;.
\label{bbbHI}\end{eqnarray} Therefore by
\eqref{d2crit}--\eqref{bbbH},
\begin{equation}\label{KNcrit} K(z_0) =
\frac{1}{\pi^{{m+2\choose 2}}\det A\det\La} \int
\left|\det(H'H'{}^*-|x|^2I)\right|e^{ -{\langle \La^{-1}(H'\oplus
x),H'\oplus x\rangle}}\,dH'\,dx\,. \end{equation} This completes
the proof of Corollary \ref{KNcrit2}.\qed

\subsection{Alternate viewpoint}\label{s-alt}

In this section, we give a different viewpoint to the proof of
Theorem \ref{KNcrit1}  that seems closer to the discussions in
\cite{Doug}.

 Let us first
consider the simpler case of critical points $\nabla f = 0$ with
respect to the usual flat Euclidean gradient  of real valued
functions  $f: \R^n \to \R$. The delta function on the critical
set is then given by
$$C_f= \sum_{x: df(x) = 0} \delta_{x},$$
where $\delta_x$ denotes the point mass at the point $x$, i.e.
$\langle \delta_x, \phi \rangle = \phi(x)$ for test functions
$\phi$.  The measure $C_f$ is closely related to the pull back
under $\nabla f: \R^n \to \R^n$ of the delta function $\delta_0$
at zero in $\R^n$. In general, let $F: \R^n \to \R^n$ be a smooth
map all of whose zeros are non-degenerate in the sense that $\det
DF_{x} \not = 0$ whenever $F(x) = 0$. Then,
\begin{equation} \label{CAPF} F^* \delta_0 = \sum_{x: F(x) =
0} \frac{\delta_{x}}{|\det DF_{x}|}. \end{equation} If $F = \nabla
f$ and $f$ has only non-degenerate critical points, this becomes
\begin{equation} \label{REALPB} (\nabla f)^* \delta_0 = \sum_{x:
df(x) = 0} \frac{\delta_{x}}{|\det D \nabla f(x)|},
\end{equation} where $D \nabla f$ is the derivative of
the map $\nabla f$. The measures $C_f$ and $\nabla f^*\delta_0$
are related by
\begin{equation} \label{REALCF} C_f = |\det D \nabla  f(x)| (\nabla
f)^*\delta_0.\end{equation}

We now generalize to the local analogue of the case of concern in
this paper, where $f: \C^m \to \C$ is holomorphic and $\nabla$ is
a smooth connection of type $(1,0)$, i.e. has the form  $\nabla f
=
\partial f -f \partial K$ for $f$ holomorphic. As before, we define
$$C_f = \sum_{z: \nabla f(z) = 0} \delta_{z}, \;\; z \in \C^m.
$$ Relative to the global  basis $dz_j$ of $(1, 0)$ forms on $\C^m$, we may
express  $\nabla f$ as the smooth map
 \begin{equation} \label{GradM} \nabla f = (\nabla_{\frac{\partial }{\partial
z_1}} f, \dots, \nabla_{\frac{\partial }{\partial z_m}}f):
 \C^m \to \C^m, \end{equation}
where $\nabla_{\frac{\partial }{\partial z_j}} f = \frac{\partial
}{\partial z_j} f - \frac{\partial K}{\partial z_j}  f. $

Since $\nabla f$ is not holomorphic, its derivative $D^v \nabla
f(z)$ is not a complex linear map on the complex tangent space
$T_z\C^m \sim \C^m$, but rather is a linear map of the real
tangent space, a real $2m$-dimensional vector space.
 At a
critical point $z_0$, we express  the derivative  $D^v \nabla f
(z): T_z \R^{2m} \otimes \C \to T_{0} \R^{2m} \otimes \C$ in terms
of the real basis $\frac{\partial }{\partial z_j}, \frac{\partial
}{\partial \bar{z}_k}$ of each complexified real tangent space:
\begin{equation}\label{Dv}  D^v \nabla f =
\left( \begin{array}{ll} D_{\frac{\partial }{\partial
z_k}}\nabla_{\frac{\partial } {\partial z_j}} f &
D_{\frac{\partial }{\partial \bar{z}_k}} \nabla_{\frac{\partial
}{\partial z_j}} {f}
\\ & \\
D_{\frac{\partial }{\partial {z}_k}}\nabla_{\frac{\partial
}{\partial \bar{z}_j}}\bar f  & D_{\frac{\partial }{\partial
\bar{z}_k}}\nabla_{\frac{\partial }{\partial \bar{z}_j}} \bar{f}
\end{array} \right). \end{equation}
Since $f$ is holomorphic, the off-diagonal blocks simplify to
$\overline{ \Theta f}$ and its complex conjugate.  Thus, $D^v \nabla f$
is precisely the complex Hessian matrix $H^c$ of (\ref{Hmatrix}), and
hence it is  the
`vertical derivative'  from \S \ref{HESS1}.

We then have
\begin{equation} \label{CXPB0} (\nabla f)^*  \delta_0= \sum_{z:
df(z) = 0} \frac{ \delta_{z}}{|\det D^v \nabla  f(z)|},
\end{equation}
and therefore
$$C_f = |\det
D^v \nabla  f(z)| \sum_{z: df(z) = 0} (\nabla f)^*  \delta_{z}, $$
where $D^v \nabla f$ is given in (\ref{Dv}).

We now adapt these formulas to  holomorphic sections $s \in H^0(M,
L)$ and  Chern connections for a Hermitian metric $h$, which
reduce to the previous example in a local frame. In so doing, we
justify the invariant interpretation of the Hermitian Hessian
matrix $H^c$ from \S \ref{HESS1}.

We introduce local coordinates $z_1,\dots,z_n$ on $M$ with
Euclidean volume form $dz$ as in Theorem \ref{KNcrit1}. From an
invariant point of view, the gradient map (\ref{NABLAMAP}) is a
section of $ T^{* (1, 0)} \otimes L.$ Since the bundle is
non-trivial, the delta-function at $0$ in the previous calculation
should be interpreted as  the delta function $\delta_0$ along the
zero section of $T^{* (1, 0)} \otimes L$, that is,
$$\langle \delta_0, \psi \rangle = \int_M \psi(z, 0) dV(z) $$
where $\psi \in \ccal^\infty(T^{*1,0} \otimes L). $ In a local
frame $e_L$,  the gradient map is given by (\ref{GradM}) and the
delta-function is just $\delta_0$ on $\C^m$.   This explains why
the derivative $D^v$ defined in the local discussion is the
vertical part of the full derivative relative to the flat
connection. In the setting of line bundles, (\ref{CXPB0}) becomes
\begin{equation} \label{CXPB} (\nabla s)^*  \delta_0 = \sum_{z:
\nabla s(z) = 0} \frac{\delta_z}{|\det D^v \nabla  s(z)|},
\end{equation}
where \begin{equation} D^v \nabla s: T_z M \to  T^{* 1, 0}_z
\otimes L_z
\end{equation}
is the vertical part of the derivative of (\ref{NABLAMAP}). Taking
the vertical part requires a connection on   $T^{* 1, 0} \otimes
L$, which we take to be the flat connection on our coordinate
neighborhood (which we can do, since $\k_{\scal,\ga,\nabla}$ is
independent of the connection on $M$). As observed in
\S \ref{HESS1}, at a critical point, $D^v \nabla s (z)$ is
independent of the choice of connection. The determinant is taken
relative to the local Euclidean metric on $M$ and the metric $h$
on $L_z.$

It follows that
\begin{equation} \label{ECS} \k_{\scal,\ga,\nabla}\,dz =
{\bf E} C_s^\nabla = {\bf E}\big(
 |\det D^v \nabla s|\, \nabla s^* \delta_0\big). \end{equation}
In calculating $\k_{\scal,\ga,\nabla}(z)$ we may fix $z$ and
introduce a local adapted frame $e_L$ at $z$. Again writing
$s=fe_L$, we have
\begin{equation}\begin{array}{rcl}
\k_{\scal,\ga,\nabla} (z)  &  = &\displaystyle  \int_{\C^m}
\int_{\C^{d}} |\det \sum a_j D^v \nabla  f_j(z)| e^{i \Re\langle t
, \sum_j a_j \nabla f_j(z) \rangle} e^{- |a|^2} dadt
\end{array}
\end{equation}
In local coordinates, $D^v \nabla f(z)$ is    the matrix $H^c$ of
(\ref{Hmatrix}).

 We now calculate the  $da$ integral by making a change of variables. We
consider the real linear  map $\jcal:=J^1_z$ of (\ref{JCAL}), which is
locally written as
\begin{equation} \jcal (a) = (\xi, H) := ( \sum_j a_j \nabla f_j(z), \sum
a_j D^v \nabla f_j(z)). \end{equation} As mentioned above, for a
positive line bundle the $H$ matrix depends only on a complex $m
\times m$ symmetric matrix (the holomorphic Hessian) and a complex
scalar (which when multiplied by $\Theta$, gives the mixed Hessian), so we
may regard
$\jcal$ as a map from $a \in \C^{d}$ into $(\xi, H) \in \C^m \times \sym(m,
\C)
\times \C$ of dimension ${m+2\choose 2}$. Since the integrand is a
function only of $\xi, H$ we may push forward the measure
$e^{-|a|^2/2} da$ under $J^1_z $ to obtain
\begin{equation}\label{REWRITE}
\int_{\C^m} \int_{\C^m \times Sym(m, \C) \times \C} |\det H|
e^{i\Re \langle t , \xi \rangle}  J(\xi, H)  d\xi dH dt,
\end{equation}
where
$$\jcal_* e^{-|a|^2} da =  J(\xi, H)  d\xi
dH, \;\; \mbox{i.e.}\; J(\xi, H)  = \int_{\jcal^{-1}(\xi, H)}
e^{-|a|^2/2} d\dot{a}
$$ where $d \dot{a}$ is the surface Lebesgue measure on the subspace $\jcal^{-1}(\xi, H)$ .
 Evaluating the $dt$ integral we obtain
\begin{equation} \mbox{(\ref{REWRITE}) } =
 \int_{\sym(m, \C) \times \C} |\det H|   J(0, H)   dH
 \end{equation}

To complete the proof, we need to evaluate  $ J(\xi, H)$. We claim
that \begin{equation}\label{J}  J(0, H)  = \frac{1}{\det A \;\det
\Lambda} e^{- \langle  \Lambda^{-1} H, H \rangle } \end{equation}
in our previous notation. This follows from general principles on
pushing forward complex Gaussians under complex linear maps $F:
\C^d \to \C^n$, whereby
$$F_* e^{-|a|^2} da = \gamma_{F F^*}, $$ i.e.
\begin{equation} \label{JCALJCAL}  J (\xi, H)   = \frac{1}{\det \jcal
\jcal^*} e^{ \langle [\jcal  \jcal^*]^{-1} (H, \xi), (H, \xi) \rangle }.
\end{equation} From
$$\langle \jcal(a), ( \xi,H) \rangle =  \{\langle \sum_j a_j\nabla
f_j(z), \xi \rangle +  \langle D^v \nabla  f_j(z), H \rangle \}$$
we see that $\jcal ^*: \C^m\times\sym(m, \C) \times \C \to \C^d$ is the
map
$$\jcal ^*(H, \xi) = (\{\langle \nabla
f_j(z), \xi \rangle + \langle D^v \nabla  f_j(z), H \rangle
\}_{j = 1}^{d}. $$ Hence $\jcal \jcal^*: \sym(m, \C) \times \C^m
\to \sym(m, \C) \times \C^m$ is the map with  block matrix form
$$\jcal \jcal^*(H, \xi)  = \left( \begin{array}{ll} A (\xi, H) & B(\xi, H) \\ & \\ B^*(\xi, H) &  C(\xi, H) \end{array}
\right)$$ where
$$\left\{ \begin{array}{l} A(\xi) = \sum_j \langle \nabla f_j(z), \xi
\rangle \nabla f_j(z), \\ \\
B(\xi, H) =  \sum_j \langle D^v \nabla  f_j(z), H \rangle \nabla
f_j(z) \oplus \langle \nabla f_j(z), \xi \rangle  D^v \nabla
f_j(z) , \\ \\
C(H)  = \sum_j \langle D^v \nabla f_j(z), H \rangle  D^v \nabla
f_j(z).
\end{array}\right.
$$
Summing in $j$ we observe that
$$A = \sum_j \langle \nabla f_j(z), \xi
\rangle \nabla f_j(z) = \nabla_z \nabla_{\bar{w}} \Pi(z, w)|_{z =
w},\;\;\; B = (\nabla_z \nabla_{\bar{w}}^2 \Pi(z,w)|_{z = w}, \;
\nabla_z \Pi(z,w)|_{z = w} ),
$$ and \begin{equation} C = \left(\begin{array}{ll} \nabla_z^2
\nabla_{\bar{w}}^2 \Pi (z, w)|_{z = w} &  \nabla_z^2\Pi(z,w)|_{z =
w}
\\ &
\\ \nabla_{\bar{w}}^2 \Pi(z, w)|_{z = w} & \Pi (z,z)
\end{array} \right).
\end{equation}
Here, $|_{z = w}$   is shorthand for $\rho^{\diag}_{e_L(z)}$ (see
(\ref{rhodiag}).

To complete the proof of Theorem  \ref{KNcrit1} we observe that
when we set $\xi = 0$, the quadratic form $\langle (\jcal
\jcal^*)^{-1}(0, H), (0, H) \rangle$ equals $\langle
\Lambda^{-1}(0, H), (0, H) \rangle$, where $\Lambda$ is given by
\eqref{Lambda}. \qed

\section{Exact formulas for Riemann surfaces}

In this section, we derive exact formulas for the density of
critical points on a Riemann surface with respect to any Hermitian
line bundle, positive or not.

\subsection{Density of critical points  on a Riemann
surface: Proof of Theorem \ref{crit1} }\label{SIMPLIFY1}

Let $(L, h) \to M$ be a  Hermitian
line
bundle on a Riemann surface $M$ with area
form $dV$, and let $\scal$ be a finite-dimensional subspace of
$H^0(M,L)$ with the 2-jet spanning property, as in Theorem
\ref{crit1}. Let $r=\frac
i2
\Theta_h/dV$, and let $\mu_1, \mu_2$ be the eigenvalues of $\La
Q_r$, where
$$Q =\begin{pmatrix} 1 & 0 \\ & \\ 0 & -r^2
\end{pmatrix}\;.$$
 We
observe that  $\mu_1, \mu_2$ have opposite signs since $\det \La
Q_r = - r^2\det \Lambda < 0$. Let $\mu_2 < 0 < \mu_1$.

From the 1-dimensional case of Theorem \ref{KNcrit1}, we have
\begin{equation}\label{dim1}\kcal^\crit_{ \scal,h,V}  =
\frac{1}{\pi^3 A \det\La} \int_{\C^2}
 \big|\,|H'|^2
-r^2| x|^2\,\big|e^{ -{\langle
\La^{-1}(H',x),(H',x) \rangle}}\,dH'\,dx\,.\end{equation}

 Writing $H=\big(\,H'\ \ x\,\big)$, we then have:
\begin{eqnarray*}\kcal^\crit_{ \scal,h,V}&=&
\frac{1}{\pi^3 A\det\La} \int_{\C^2} \big|HQ_rH^*\big|\exp\left(
-H\La\inv H^*\right)\,dH
\\&=&
\frac{1}{\pi^3 A} \int_{\C^2} \big|H\La^\half Q_r\La^\half
H^*\big|\exp\left( -H H^*\right)\,dH\,.\end{eqnarray*}
We diagonalize $\La^\half Q_r\La^\half $, which has the same
eigenvalues $\mu_1,\ \mu_2$ as
$\La Q_r$, to obtain
\begin{eqnarray*}\kcal^\crit_{ \scal,h,V}&=&
\frac{1}{\pi^3A}\int_{\C^2}\big|\,\mu_1|a|^2 +\mu_2|b|^2\,\big|\; e^{-
|a|^2 -|b|^2}\;da\,db \\
&=& \frac{1}{\pi A}\int_0^{+\infty}\! \int_0^{+\infty}
|\mu_1 u+\mu_2v|\;e^{- u -v}du\,dv\\
&=& \frac{1}{\pi A\mu_1|\mu_2|}\int_0^{+\infty}\! \int_0^{+\infty}
| u-v|\;e^{- \mu_1\inv u -|\mu_2|\inv v}du\,dv\\
&=&
 \frac{1}{\pi A\mu_1|\mu_2|}\int_{-\infty}^{+\infty}\int_{\max\{w,0\}}^{+\infty}
|w|\;\exp(|\mu_2|\inv  w)
\exp\left[-(\mu_1\inv +|\mu_2|\inv)
u\right]\;du\;dw\\
&=&  \frac{1}{\pi A\mu_1|\mu_2|}\;\mbox{(I)}\ +\  \frac{1}{\pi
A\mu_1|\mu_2|}\;\mbox{(II)}\,,
\end{eqnarray*}
where
\begin{equation*} \mbox{(I)} =
\int_{0}^{+\infty}\int_{w}^{+\infty} w\;\exp(|\mu_2|\inv  w)
\exp\left[-(\mu_1\inv +|\mu_2|\inv)
u\right]\;du\;dw = \frac{\mu_1^2}{\mu_1\inv +|\mu_2|\inv}\;,
\end{equation*} and
\begin{equation*} \mbox{(II)} = \frac{1}{2\pi}\int_{-\infty}^{0}\int_{0}^{+\infty}
(-w)\;\exp(|\mu_2|\inv  w)
\exp\left[-(\mu_1\inv +|\mu_2|\inv)
u\right]\;du\;dw = \frac{\mu_2^2}{\mu_1\inv +|\mu_2|\inv}\;.\end{equation*}
This yields the desired formula.\qed

\subsection{Index density: Proof of Corollary \ref{indexdensity}}

Critical points of a section $s$ are (almost surely) of index $\pm
1$. The above proof shows that the expected density of critical
points of index $1$ is given by $$\kcal^\crit_+= \frac{1}{\pi
A\mu_1|\mu_2|}\;\mbox{(I)} = \frac{1}{\pi A}\; \frac{\mu_1^2
}{| \mu_1| + | \mu_2|} \;,$$ while the expected density of
critical points of index $-1$ is
$$\kcal^\crit_-= \frac{1}{\pi A\mu_1|\mu_2|}\;\mbox{(II)} =
\frac{1}{\pi A}\; \frac{\mu_2^2  }{| \mu_1| + | \mu_2|} \;.$$
Hence, the index density is given by
\begin{equation}\label{index1}\kcal^\crit_{\rm index}: =
\kcal^\crit_+ - \kcal^\crit_- =
\frac{\mu_1+\mu_2}{\pi A} = \frac{Tr[\La Q_r]}{\pi A}
\;.\end{equation} (Of course, \eqref{index1} can also be obtained
directly from \eqref{dim1}  as an elementary second moment
calculation.)

The critical points of $s$ of index $1$ are the saddle points of
$|s|^2_h$, while those of index $-1$ are local maxima of $|s|^2_h$
in the case where $L$ is positive, and are local minima of
$|s|^2_h$ if $L$ is negative. (The length $|s|$ cannot have positive
local minima if $L$ is positive, or maxima if $L$ is negative.)

\subsection{\label{pfcrit1}Alternate proof of Theorem \ref{crit1}} For
simplicity, we assume that $r=1$.  From \eqref{dim1}, we obtain
\begin{eqnarray*}\kcal^\crit_{ \scal,h,V}  &=&
\frac{1}{\pi^3 A \det\La} \int_{\C^2}
 \big|\,|H'|^2
-| x|^2\,\big|\,e^{ -{\langle
\La^{-1}(H',x),(H',x) \rangle}}\,dH'\,dx
\\&=&
 \frac 1 {2\pi^4A\, \det\La}\int_{
\C ^2} \int_{\R} \int_{\R}\left|\lambda-|x|^2 \right| e^{i  \xi
(\lambda- |H'|^2)} e^{ -{\langle
\La^{-1}(H',x),(H',x) \rangle}}d\xi\,d\la  \, dH'\, dx\,.
\end{eqnarray*}
Indeed, the $\xi$ integral gives the $\de$-function at $\la=|H'|^2$, and the
$\la$ integral then gives  \eqref{dim1}.

We change variables to $\lambda' = \lambda - |x|^2$  to get (dropping
the primes)
\begin{equation}
 \frac{1}{2\pi^4 A\det\La} \int_{ \C ^2}
\int_{\R} \int_{\R}\left|\lambda\right| e^{i  \xi( \lambda + |x|^2
-  |H|^2 )}  \exp\left( -{\langle \La^{-1}(H, x),(H, x)
\rangle}\right)  d\xi\,d\la  \, dH\,dx.
\end{equation}

We now do the complex Gaussian  $dH dx $ integral on $\C^2$. The
quadratic form is
$$ i  \langle \xi, |x|^2 - |H|^2 \rangle - \langle
\La^{-1}(H, x),(H, x) \rangle = - \left\langle (\La^{-1} + i \xi
Q) (H, x),(H, x) \right\rangle.
$$
The result  is
\begin{equation} \label{exact1}\kcal^\crit _{\scal,h,V}(z) =
 \frac{1}{2\pi^2 A}
\int_{\R} \int_{\R}\left|\lambda\right| e^{i \xi \lambda }
\frac{1}{\det(I + i \xi \Lambda Q)}\,d\xi \,d\lambda\,.
\end{equation}
Thus, in dimension one, $\kcal^\crit _{\scal,h,V}(z)$  depends only  on the
eigenvalues $\mu_1, \mu_2$ of  $\Lambda Q$.

We first consider the  $d \xi$ integral,
$$\ical(\lambda) = \int_{\R}\frac{ e^{i \xi \lambda
} }{\det(I + i \xi \Lambda Q)} d \xi = \int_{\R}\frac{ e^{i \xi
\lambda } }{(1 + i \xi \mu_1) (1 + i \xi \mu_2)} d \xi=
\frac{-1}{\mu_1 \mu_2 } \int_{\R} \frac{e^{i \xi   \lambda }}{( \xi-i
\mu_1^{-1} ) ( \xi-i \mu_2^{-1} )} d \xi. $$ We separately
treat the cases $\lambda > 0, \lambda < 0$.
\medskip

\noindent{\bf (i) $\lambda > 0$:} In this case, we pick up the
residue at the pole $i/\mu_1$ in the upper half plane:
\begin{equation} \ical_+(\lambda) = \frac{-2\pi i}{\mu_1 \mu_2}
\mbox{Res}_{i/\mu_1} \left[\frac{e^{i \xi  \lambda }}{( \xi-i
\mu_1^{-1} ) ( \xi-i \mu_2^{-1} )} \right]
 = \frac{-2\pi i}{\mu_1 \mu_2} \frac{e^{-\la/\mu_1}}
 {(i \mu_1^{-1}- i \mu_2^{-1})}
 = \frac{2\pi\, e^{-\la/\mu_1}}{\mu_1 - \mu_2} \end{equation}
 \medskip

 \noindent{\bf (ii)  $\lambda < 0$:} In this case we pick up the residue at $i/\mu_2$:
\begin{equation} \ical_-(\lambda) = \frac{2\pi i}{\mu_1 \mu_2} \mbox{Res}_{i/\mu_2}
\left[\frac{e^{i \xi  \lambda }}{( \xi-i
\mu_1^{-1} ) ( \xi-i \mu_2^{-1} )} \right] =
\frac{2\pi\, e^{\la/|\mu_2|}}{\mu_1 - \mu_2} \end{equation}

 \medskip

 To complete the calculation, we need to evaluate
 \begin{eqnarray}   \int_{-\infty}^0 (-\lambda)
\ical_-(\lambda)d\lambda  + \int_0^{\infty}
 \lambda \ical_+(\lambda) d\lambda \hspace{-1.5in}\nonumber\\& = &
  \frac{2\pi}{\mu_1 - \mu_2}\left( \int_{-\infty}^0 (-\lambda) e^{
 \frac{\lambda}{|\mu_2|}}d\lambda +
 \int_0^{\infty}  \lambda e^{-
\frac{\lambda}{\mu_1}} d\lambda\right)\nonumber\\& = &
2\pi\,\frac{\mu_2^2 + \mu_1^2 }{|\mu_1| + |\mu_2|}\;.
\label{exact2}
\end{eqnarray}  The desired formula follows from
\eqref{exact1} and \eqref{exact2}
\qed

\medskip We shall use this approach in \cite{DSZ}  for our
computation of densities in higher dimensions.

\subsection{Exact formula for $\CP^1$:
Proof of Corollary \ref{exactP1}.}\label{s-exactP1}

Since the critical point  density with respect to the Fubini-Study
metric is $SU(2)$ invariant and hence constant, it suffices to
compute it  at the point $(z_0:0)\in \CP^1$, using  the local
coordinate $z=z_1/z_0$ and the local frame $e_N$  for $\ocal(N)$
corresponding to the homogeneous polynomial $z_0^N$ on $\C^2$. We
recall that the \szego kernel is given by
$$ \Pi_{H^0(\CP^1,\ocal(N))}(z,w) = \frac {N+1}{\pi}
(1+z \bar w)^N e_N(z)\otimes \overline{e_N(w)}\;.$$ (See, for
example, \cite[\S 1.3]{SZ}.) Since our formula is invariant when
the \szego kernel is multiplied by a constant, we can replace the
above by the {\it normalized \szego kernel}  \begin{equation}
\label{NSK} \wt \Pi_N(z,w):=(1+z \bar w)^N \end{equation} in our
computation.

Since
$$K(z)=-\log|e_N(z)|^2_\FS=N\log (1+|z|^2)\;,$$ we have
$$K(0) = \frac {\d K}{\d z} (0)
=\frac {\d^2K}{\d z^2} (0) = 0\;; $$ i.e., $e_N$ is an adapted
frame at $z=0$. Hence when computing the (normalized) matrices
$B_N,\ C_N$ for $H^0(\CP^1,\ocal(N))$, we can take the usual
derivatives of $\wt\Pi_N$. Indeed, we have
\begin{eqnarray*} \frac{\d\wt\Pi_N}{\d{z}} &=&
N(1+z \bar w)^{N-1}\bar w\;,\\
\frac{\d^2\wt\Pi_N}{\d{z}\d\bar w} &=& N(1+z\bar w)^{N-1} +
N(N-1)(1+z\bar w)^{N-2}z\bar w\;.
\end{eqnarray*}
It follows that
\begin{equation}\label{ABC}A_N(0)=(\; N\; )\;,
\quad B_N(0)= (\; 0 \ \ 0\; )\;,\quad \La_N(0)=C_N(0)=
\begin{pmatrix} 2N(N-1) &0\\ 0 & 1\end{pmatrix}\;.
\end{equation}

We now apply Corollary \ref{indexdensity}. Since $r:=\frac i2
\Theta_h/dV=N$, where $h$ is the Fubini-Study metric on $\ocal(N)$
and $dV=\om_\FS$, the eigenvalues of $\La_N(0)Q_r$ are given by:
$$\mu_1=2N(N-1),\quad \mu_2=-N^2\;.$$

Suppose that $N\ge 2$ so that $\ocal(N)\to\CP^1$ has the 2-jet
spanning property. Theorem \ref{crit1} then yields
\begin{equation} \label {KP1}{\mathcal  K}^\crit_+= \frac 1\pi\;
\frac{4(N-1)^2}{3N-2}\;,\quad {\mathcal K}^\crit_- =\frac 1\pi\;
\frac{N^2}{3N-2}\;.\end{equation} Since $\kcal_\pm^\crit$ is
constant by invariance of the metric and connection, and
$\vol(\CP^1)=\pi$, the desired formulas follow from \eqref{KP1}.

If $N=1$, then every section has exactly 1 critical point (of
index 1), so the formula holds in this case too. \qed

\begin{rem} For $N=2$, it also turns out that almost all sections have
exactly 2 critical points---one each of index $+1$ and $-1$.  To see
this, we first note that Theorem
\ref{exactP1} says that the expected number of critical points in this case
is 2. Since
$\chi(\ocal(2)\otimes T^{*1,0})=c_1(\ocal(2)\otimes T^{*1,0})=0$, the number
of critical points of index 1 equals the number of critical points of index
$-1$.  Suppose that
$$s=(a+bz+cz^2)\,e^{\otimes 2}\;.$$
The critical point equation is:
$$(2a+b) +(2b+2c)z+b|z|^2=0\;.$$
By B\'ezout's Theorem on $\R^2$, there are at most 4 critical
points. Hence there are only two possibilities: (i) 2 critical
points of index 1 and 2 of index $-1$; (ii) 1 critical point each
of index 1 and of $-1$.  Since the average number of critical
points is 2, case (ii) almost always occurs.

However, for $N\ge 3$, one easily checks that the expected number of
critical points, $\frac{5N^2-8N+4}{3N-2}$, is not an integer and hence the
sections in
$H^0(\CP^1,\ocal(N))$  cannot all have the same number of critical points.
\end{rem}

\subsubsection{Metric dependence of the number of critical
points}\label{s-dependence}

The expected number of critical points $\ncal^\crit_{\scal,h}=
\int _M\kcal^\crit_{\scal,h}dV$ of a section $s$ of $H^0(M,L)$
(with the Hermitian Gaussian measure) depends on the metric $h$ on
$L$. This is true even for the case where $L=\ocal(1)\to
\CP^1$ is the hyperplane section bundle over the projective line.
To illustrate this dependence, we let $z=z_1/z_0$ denote the
coordinate in the affine chart $\C\subset\CP^1$, and let $e_L=z_0
\in \mbox{Hom}(\C^2,\C)\approx H^0(\CP^1,L)$; then $e_L$ is a
local frame over $\C$. If we give $L$ the standard Fubini-Study
metric $h(e_L,e_L)=(1+|z|^2)\inv$, then $e_L$ has a critical point
(maximum of $h(e_L,e_L)$) at $0$ and no others. Furthermore, every
section of $H^0(\CP^1,L)$ has exactly 1 critical point, so the
expected number of critical points equals 1.  Now let $p(z)$ be a
polynomial of degree $k>1$ with distinct roots, and  consider the
metric
$$\tilde h= \tilde h_0^{1-\ep}h^\ep\,, \qquad \tilde h_0 =
(p^*h)^{\frac 1k} \,,$$ where $\ep>0$. The metric $\tilde h$ has
positive curvature (while the curvature of $\tilde h_0$ is
semi-positive). Since the critical points of a section $s$
coincide with the critical points of the function $\log|s|$,  it
suffices to consider
$$\log |e_L|_{\tilde h_0} = - \frac 1 k \log  (1+|p(z)|^2)\;.$$
This function has maxima at the $k$ roots of $p$ (and for generic
$p$, has $k-1$ saddle points, by Morse theory) and hence  has
$2k-1$ critical points on $\C$.  Therefore, the section $e_L=z_0$
has $2k-1$ critical points, and hence all nearby sections
$z_0+\delta z_1$ also have $2k-1$ critical points. As every
section has at least one critical point, the expected number of
critical points is greater than 1 for the metric $\tilde h$.

On the other hand, in \cite{DSZ} we  show that for any metric with positive
curvature on
$\CP^1$ (or more generally on any compact \kahler manifold), the
expected number of critical points of $H^0(\CP^1,\ocal(N))$ has an
asymptotic expansion in $N\inv$, where the first two terms are
independent of the metric; see \eqref{critRS1}--\eqref{critRS2}.

\section{Morse index density: Proof of Theorem \ref{Morse}}\label{s-morse} We recall
that the critical points of $\log|s|^2$ coincide with the critical
points of $\nabla s$ and that they have Morse index $\ge m$.
(Since almost all sections have only nondegenerate critical
points, we make this assumption throughout.)  Theorem \ref{Morse}
is an immediate consequence of Lemma \ref{index} below and the
proof of Theorem \ref{KNcrit1} and Corollary \ref{KNcrit2}.

Recall that  the Hermitian metric on $T^{1,0}_M$ is given by the
curvature
$$\Theta=-\nabla'\nabla'' \log h\in
T^{1,0\dual}_M \otimes T^{0,1\dual}_M =T^{1,0\dual}_M \otimes
\overline{T^{1,0\dual}_M}\;.
$$ We let $\Theta^* \in
T^{1,0}_M \otimes \overline{T^{1,0}_M} $ denote the dual metric on
$T^{1,0\dual}_M$.

\begin{lem}\label{index} Let $(L,h)\to M$ be a positive holomorphic line
bundle, and let $z_0\in M$ be a nondegenerate critical point of
$s\in H^0(M,L)$. Then the Morse index of $\log|s|_{h}^2$ at $z_0$
equals $m+\mbox{\rm index}_{z_0}(S\Theta^*\overline S -\Theta)$,
where
$$S=\nabla'\nabla'\log| s|^2_{h}\in T^{1,0\dual}_M\otimes
T^{1,0\dual}_M.$$  Hence at a critical point, the topological
index of $s$ is $(-1)^{m+n}$, where $n$ is the Morse index of
$\log|s|_h$.
\end{lem}

\begin{proof}  Let $z_0$ be a nondegenerate critical point of
$\log|s|^2$  at $z_0$. Let $z_j=x_j+iy_j,\ 1\le j\le m$, be normal
coordinates at $z_0$. Note that $\Theta_{z_0}=\sum_{j=1}^m dz_j
\otimes d\bar z_j$.  Thus
$$\mbox{\rm index}_{z_0}(S\Theta^*\overline S
-\Theta)=\mbox{\rm index}(SS^* -I)\;,$$ where $S$ now denotes the
symmetric matrix $$\big(S_{jk}\big)= \left( \frac
{\d^2}{\d z_j \d
 z_k}  \log|s|_{h}^2\right)_{z_0}\;.$$

Conjugating the Hessian matrix
$$\begin{pmatrix}\ \left(\frac {\d^2 }{\d x_j \d x_k}\log|s|_h^2
\right) &\left(\frac {\d^2}{\d x_j \d y_k} \log|s|_h^2
\right)\\[8pt]
\left(\frac {\d^2}{\d y_j \d x_k} \log|s|_h^2 \right)
&\left(\frac {\d^2}{\d y_j \d y_k} \log|s|_h^2 \right) \
\end{pmatrix}$$ with the unitary matrix $$\frac 1 {\sqrt 2}
\begin{pmatrix} I & iI\\ iI &I\end{pmatrix}\ ,$$ we get
\begin{equation}\label{matrix}2\begin{pmatrix}\
\left(\frac {\d^2}{\d z_j \d \bar z_k}  \log|s|_h^2\right)
&\left(i\frac {\d^2}{\d z_j \d z_k} \log|s|_h^2
\right)\\[8pt]
\left(-i\frac {\d^2}{\d \bar z_j \d \bar z_k} \log|s|_h^2
\right) &\left(\frac {\d^2}{\d z_j \d \bar z_k}
\log|s|_h^2\right) \
\end{pmatrix}\ .\end{equation}
Write $s=f\,e_L$, where $e_L$ is a local frame for $L$.
Since
$$\log |s|^2_h= \log |f|^2 + \log h$$ and $\{z_j\}$ are normal
coordinates, we have
$$\frac {\d^2}{\d z_j \d \bar z_k}
\log|s|_h^2 =-\delta_j^k \;.$$

Thus,  \eqref{matrix} becomes:
\begin{equation}\label{matrix*}2\begin{pmatrix}
-I & iS\\ -iS^* & -I
\end{pmatrix}\ .\end{equation}
Let $\lambda$ be an eigenvalue of $$\wh S:=\begin{pmatrix} -I &
iS\\ -iS^* & -I
\end{pmatrix}\ ,$$  Then
$$0=\det \begin{pmatrix}
(-1-\la)I & iS\\ -iS^* & (-1-\la)I
\end{pmatrix}\ = \det \big[ (1+\la)^2I-SS^*\big].$$
Therefore,  $\mu:=(1+\la)^2$ is an eigenvalue of $SS^*$. On the
other hand, each eigenvalue $\mu$ of $SS^*$ yields the pair of
eigenvalues $-1\pm \sqrt\mu$ of $\wh S$. Hence the number of
negative eigenvalues of the Hessian of $\log|s|_h^2$ equals
$m$ plus the number of eigenvalues of $SS^*$ that are less than 1.
\end{proof}

\begin{rem}
If the line bundle $(L,h)\to M$ instead has negative curvature, then
the Morse index of $\log|s|_h^2$ at a critical point $z_0\in M$
equals\/ $\mbox{\rm index}_{z_0}(S\Theta^*\overline S -\Theta)$. To see
this, we choose normal coordinates at $z_0$ with respect to the metric
$-\Theta$.  Then,
$$\mbox{\rm index}_{z_0}(S\Theta^*\overline S
-\Theta)=\mbox{\rm index}(I-SS^*)\;.$$ This time, the Morse index
of $\log|s|_h^2$ is the number of negative eigenvalues of
$$\wh S:=\begin{pmatrix}
I & iS\\ -iS^* & I
\end{pmatrix} .$$  Thus, each eigenvalue $\mu$ of
$SS^*$ corresponds to the pair of eigenvalues $1\pm \sqrt \mu$ of
$\wh S$, and the conclusion follows as above.\end{rem}

\end{document}